\documentclass[a4paper,11pt]{amsart}
\usepackage{amsmath, amssymb, amsfonts, amsthm, enumerate, mathrsfs, array, bm}

\newtheorem{thm}{Theorem}
\newtheorem{cor}[thm]{Corollary}
\newtheorem{lem}[thm]{Lemma}

\newtheorem{prop}[thm]{Proposition}

\newtheorem{mainthm}[thm]{Main Theorem}

\newtheorem{rem}{Remark}

\newtheorem{hyp}{Hypothesis}[thm]

\newcommand{\Q}{\mathbb Q}
\newcommand{\R}{\mathbb R}
\newcommand{\C}{\mathbb C}
\newcommand{\Z}{\mathbb Z}
\newcommand{\A}{\mathbb A}
\newcommand{\G}{\mathbb G}

\newcommand{\T}{\mathbb T}
\newcommand{\OS}{\mathscr O}
\newcommand{\SGU}{S(G,U)}
\newcommand{\tildeSGU}{S(\tilde{G}, \tilde{U})}
\newcommand{\tildeSTU}{S(\tilde{T}, \nu(\tilde{U}))}
\newcommand{\STU}{S(T, \nu(U))}
\newcommand{\GLn}{GL_n}
\newcommand{\Gint}{\mathbf{G^{int}}}
\newcommand{\Hint}{\mathbf{H^{int}}}

\numberwithin{equation}{section}

\begin{document}

\title{Level-raising for automorphic forms on $GL_n$}

\author{Aditya Karnataki}
\address{Beijing International Center for Mathematical Research \\ Peking University \\ 5 Yiheyuan Road \\ Beijing \\ People's Republic of China 100871}

\email{adityack@bicmr.pku.edu.cn}
\thanks{I am deeply grateful to Jack Thorne for suggesting the project underlying this work and for his constant support. I am thankful to the anonymous referee for their careful reading of this article and suggestions. Additionally, It is a pleasure to thank David Rohrlich, Ben Howard, Dipendra Prasad, C. S. Rajan, Liang Xiao for their useful comments. I am thankful to Yakov Varshavsky for directing my attention to results in his paper that simplified the proof immensely. Parts of this work were produced when the author was a visiting fellow at Tata Institute of Fundamental Research, Mumbai, and I would like to thank the institute for its hospitality.}


\date{}

\dedicatory{}

\begin{abstract}
Let $E$ be a CM number field and $F$ its maximal real subfield. We prove a level-raising result for regular algebraic conjugate self-dual automorphic representations of $\GLn(\A_E)$. This generalizes previously known results of Thorne \cite{Thorne} by removing certain hypotheses occurring in that work. In particular, the level-raising prime $p$ is allowed to be unramified as opposed to inert in $F$, the field $E/F$ is not assumed to be everywhere unramified, and the field $E$ is allowed to be a general CM field. 
\end{abstract}

\maketitle

\section{Introduction}
\label{sec:introduction}

Let $f$ be a normalized elliptic eigenform for the group $\Gamma_{0}(N)$ of weight two. Then, %\begin{equation*} f(z) = \sum_{n=1}^{\infty} a_n(f) e^{2 \pi i z} \end{equation*} with $a_1(f) = 1$, and the Fourier coefficients $a_n(f)$ are eigenvalues of the Hecke operators $T_n$. If $\ell$ is a prime, and $\iota : \bar{\Q}_{\ell} \cong \mathbb{C}$ is a chosen isomorphism, 
there exists an associated Galois representation \begin{equation*}
    r_{\iota}(f) : G_{\mathbb{Q}} \rightarrow GL_{2}(\bar{\Q}_{\ell}),
\end{equation*} unramified outside $N\ell$, uniquely characterised by the fact that the trace of Frobenius at a finite prime $p$ not dividing $N \ell$ equals the Fourier coefficient $a_p(f)$. One can even assume that $r_{\iota}(f)$ is integral without loss of generality by choosing an appropriate basis, and thus we may consider its reduction $\overline{r_{\iota}(f)} : G_{\Q} \rightarrow \bar{\mathbb{F}}_{\ell}$. We suppose that this is irreducible. Let $p$ be a prime away from $N \ell$. If the so-called `level-raising congruence' $a_p(f)^2 \equiv (p+1)^2 \mod{\ell}$ holds, then locally at $p$, there exists a lift of $\overline{r_{\iota}(f)}|_{G_{\Q_p}}$ to a representation $\rho : G_{\Q_p} \rightarrow GL_{2}(\bar{\mathbb{Z}}_{\ell})$ such that, under the local Langlands correspondence, $\rho \otimes_{\bar{\mathbb{Z}}_{\ell}} \bar{\mathbb{Q}}_{\ell}$ corresponds to an unramified twist of the Steinberg representation. This naturally leads to the idea of asking whether there exists a \textit{global} tamely ramified Galois representation congruent to $r_{\iota}(f)$ modulo $\ell$, and consequently, whether there exists a modular form with $p$ dividing level exactly once satisfying this. The level-raising congruence is necessary due to local-global compatibility obstruction. Then one can ask about the existence of an elliptic $p$-new eigenform $g$ for the group $\Gamma_{0}(Np)$ and of weight two, such that $\overline{r_{\iota}(g)} \cong \overline{r_{\iota}(f)}$.

In this classical case for elliptic modular forms, this question was first asked and answered affirmatively by Ribet \cite{Ribet}. Since the level of the modular form in question goes from $N$ to $Np$, this is a question of raising the level. Since Ribet and Hida's pioneering work, congruences between automorphic representations has been a central theme in Number Theory. Level-raising congruences in particular have played an essential role in several key developments in areas such as automorphy lifting theorems and arithmetic of elliptic curves. In recent times, establishing some key results in this direction such as \cite{AT} has allowed for proving the symmetric power functoriality for holomorphic modular forms as pursued in \cite{NT}.

One of the key inputs in \cite{Ribet} is a lemma of Ihara, which describes the kernel of a degeneracy map \begin{equation*}
    H^1(\Gamma_{0}(N), \mathbb{F}_{\ell}) \oplus H^1(\Gamma_{0}(N), \mathbb{F}_{\ell}) \rightarrow H^1(\Gamma_{0}(Np), \mathbb{F}_{\ell}).
\end{equation*} In its essentials, this is a statement about the group $GL_2$ and its arithmetic subgroups. Diamond and Taylor \cite{DT} generalized \cite{Ribet} by considering automorphic forms on quaternion algebras, i.e. inner forms of $GL_2$ over $\Q$. Ihara's lemma for groups of higher rank is less well-understood in general. Clozel, Harris, and Taylor \cite{CHT} formulated a conjecture about the structure of space of automorphic forms on definite unitary groups, that would play the role, for applications to level-raising, of Ihara's Lemma in Ribet's work. In \cite{Thorne}, Thorne establishes instances of a conjectural `Ihara lemma' of \cite{CHT} for unitary groups of rank higher than $2$, thereby simultaneously proving a level-raising result for regular algebraic, conjugate self-dual automorphic representations of $\GLn(\A_E)$, where $E$ is a CM number field.

In present article, we seek to prove new level-raising theorems in banal characteristic, for automorphic representations $\pi$ of $GL_n(\A_E)$, where $E$ is a CM number field, i.e. a totally imaginary extension of a totally real number field, by seeking to remove some of the hypotheses present in \cite{Thorne}. We suppose that $\pi$ is regular algebraic and conjugate self-dual. In this case, the existence of a Galois representation $r_{\iota}(\pi) : G_E \rightarrow GL_n(\bar{\Q}_{\ell})$ is known, and the question of level-raising can be formulated much the same way as was done for elliptic modular forms.  We first state a version of our main theorem of level-raising in this context below. This is Theorem \ref{final} in the body of the paper. 

Let $E$ be a CM field with maximal totally real subfield $F$, $[F:\Q] = d$ is even. The level-raising prime $p$ is assumed to be unramified in $F$, with a place $\nu$ above it of $F$ assumed to be split in $E$ as $\nu = \omega . \omega^{c}$. Let $n \ge 3$ be an integer, and $\ell \not = p$ a prime. We fix an isomorphism $\iota : \bar{\Q}_{\ell} \cong \C $. Let $n_1, n_2$ be positive integers with $ n = n_1 + n_2 $. Suppose that there exist $\pi_1, \pi_2$ conjugate self-dual cuspidal automorphic representations of $GL_{n_1}(\A_E)$ and $GL_{n_2}(\A_E)$ respectively such that $\pi = \pi_1 \boxplus \pi_2$ is regular algebraic. Let $\bm{\lambda} = (\lambda_{\tau})_{\tau : E \hookrightarrow \C}$ be the weight of $\pi$. Under these hypotheses, it is known that there is associated to $\pi$ a continuous semisimple Galois representation $r_{\iota}(\pi) : G_{E} \rightarrow GL_n(\bar{\Q}_{\ell})$.

\begin{mainthm}

With $\pi$ as above, suppose that the local component at $\omega$ of $\pi$, denoted as $\iota^{-1} \pi_{\omega}$ satisfies a level-raising congruence \ref{eq:LRC} modulo $\ell$ which is assumed to be a banal characteristic. Suppose further that the following hypotheses hold.

\begin{enumerate}
\item If $t_{\ell}$ is a generator of the $\ell$-th part of the tame inertia group at $\omega$, then $\overline{r_{\iota}(\pi)}(t_{\ell})$ is a unipotent matrix with exactly two Jordan blocks.

\item The weight $\bm{\lambda} = (\lambda_{\tau})_{\tau : E \hookrightarrow \C}$ of $\pi$ satisfy the following properties :
          
            \begin{itemize}

            \item For each $\tau$, and for each $0 \le i < j \le n$, we have $0 < \lambda_{\tau, i} - \lambda_{\tau, j} < \ell$.

            \item There exists an isomorphism $\iota_{p} : \overline{\Q}_p \cong \C$ such that the following inequalities hold : 

		$$ 2n + \sum_{\tau : E \hookrightarrow \C} \sum_{j = 1}^{n} (\lambda_{\tau, j} - 2 \left \lfloor{ \lambda_{\tau, n}/2 }\right \rfloor ) \le \ell,$$
		$$ 2n + \sum_{\tau : E \hookrightarrow \C} \sum_{j = 1}^{n} (2 \left \lceil{ \lambda_{\tau, 1}/2 }\right \rceil  - \lambda_{\tau, n+1-j}) \le \ell.$$

            \end{itemize}

\item If $\pi$ is ramified at a place $\gamma$ of $E$, then $\gamma$ is split over $F$.

\item $\pi$ is unramified at the primes of $E$ lying above $\ell$, and $\ell$ is unramified in $E$.

\item $\pi = \pi_1 \boxplus \pi_2$ satisfies a sign condition \ref{eq:sgn}, $n_1 \not = n_2$, and $n_1n_2$ is even.

\end{enumerate}
Then there exists a regular algebraic conjugate self-dual cuspidal automorphic representation $\Pi$ of $GL_n(\A_E)$ of weight $\bm{\lambda}$ such that there is an associated Galois representation $r_{\iota}(\Pi)$ to it satisfying $\overline{r_{\iota}(\pi)} \cong \overline{r_{\iota}(\Pi)}$ and $\Pi_{\omega}$ is an unramified twist of the Steinberg representation.

%With $\pi$ as above, suppose that $\iota^{-1} \pi_{\omega}$ satisfies the level-raising congruence $(4.2)$ of \cite{Thorne}. Then, under certain hypotheses, there exists a regular algebraic conjugate self-dual cuspidal automorphic representation $\Pi$ of $GL_n(\A_E)$ of weight $\bm{\lambda}$ such that there is an associated Galois representation $r_{\iota}(\Pi)$ to it satisfying $\overline{r_{\iota}(\pi)} \cong \overline{r_{\iota}(\Pi)}$ and $\Pi_{\omega}$ is an unramified twist of the Steinberg representation. 
\end{mainthm}

 %Thorne uses crucially the Shimura varieties arising from unitary similitude groups for establishing his result. Use of these varieties results in some additional hypotheses on the class of fields $E$ appearing in that work. Inspired by his approach, our aim is to generalize his results for a larger class of cases by removing these hypotheses, by using certain locally symmetric spaces that arise from true unitary groups (as opposed to unitary similitude groups), and proving $p$-adic uniformization results for such models. Then we apply these results to the problem of raising the level for automorphic forms.

This theorem improves \cite[Theorem 7.1]{Thorne} as mentioned earlier. There $p$ is assumed to be inert in the totally real field $F$, while we allow it to be unramified. The CM field $E$ is not assumed to contain an imaginary quadratic field. It is also not required that $E/F$ be everywhere unramified.

One of the main interests in such a level-raising theorem comes from the fact that it can often allow one to prove automorphy of Galois representations $r : G_E : GL_n(\bar{\Q}_{\ell})$ satisfying the kinds of conditions as follows -

\begin{itemize}
    \item $r$ is ordinary, and there exists a place $\omega$ of $E$ with residue characteristic different from $\ell$ at which $r$ appears to correspond to a twist of Steinberg representation.
    
    \item The residual representation $\bar{r}$ is reducible, and the Jordan-Holder factors of $\bar{r}$ are individually residually automorphic.
\end{itemize}

The conditions $E$ contains an imaginary quadratic field and $E/F$ everywhere unramified can be assumed for applications to such automorphy lifting theorems by making a series of cyclic base changes following the book of Arthur-Clozel (\cite{AC}), but the condition on $p$ being unramified should lead to improved automorphy lifting theorems. We also mention here that level-raising results of this type have also had applications to generalizations of Gross-Zagier-Kolyvagin type results, \cite{Zhang} for instance, and here it is important to be able to prove something without making an extension of the base field.

We briefly describe the main ideas of our work. Fix an integer $n > 1$. Let $E$ be an imaginary CM number field. By definition, $E$ is an imaginary quadratic extension of a totally real subfield, say $F$ with $[F: \Q] = d$. Let $\pi$ be a regular algebraic conjugate self-dual cuspidal representation of $GL_n(\A_E)$. If $I$ is a definite unitary group in $n$ variables associated to the extension $E/F$, then it is known that $\pi$ can be `transferred' to an automorphic representation $\sigma$ of $I(\A_F)$. Such a representation of $I(\A_F)$ can be described in terms of spaces of algebraic modular forms, which admit integral structures. In \cite{Thorne}, this representation is shown to further correspond to a representation $\tilde{\sigma}$ on $\tilde{I}$, which is defined as a group of unitary similitudes related to $I$. Let $\tilde{G}$ be an inner form of $\tilde{I}$ which differs from $\tilde{I}$ only at a fixed prime $p$ and $\infty$. $\tilde{G}$ has type $U(n-1, 1) \times U(n)^{d-1}$ at infinity, and it looks locally like the unit group of a division algebra at $p$. We note that inner forms obtained by such interchange of invariants appear in the case of $p$-adic uniformizations of Shimura curves, studied first by Cerednik \cite{Cerednik}. Drinfeld \cite{Drinfeld} later found a natural proof of Cerednik's result using moduli interpretations. Rapoport and Zink \cite{RZ} generalized these results for Shimura varieties wherever such a moduli interpretation was possible. In particular, their results apply to the groups of unitary similitudes $\tilde{I}$ and $\tilde{G}$. These results were used by Thorne to relate the space of algebraic automorphic forms on $\tilde{I}$, and $\tilde{\sigma}$, to the cohomology of Shimura variety arising from $\tilde{G}$. He then proves his level-raising result using the weight spectral sequence. Unfortunately, passage from $GL_n$ to $\tilde{I}$ forces some extra hypotheses on the number field $E$ at the beginning. For instance, $E$ is required to have the form $E = E_0.F$, where $E_0$ is an imaginary quadratic field. Similarly, $p$ is required to be inert in $F$. We remove in particular both of these hypotheses in present work. We note that we require $p$ to be unramified in $F$, which is an artifact of the method of using $p$-adic uniformizations.

The first new ingredient in this work is to prove $p$-adic uniformization results for some locally symmetric spaces arising from the ``true" unitary subgroup $G \subset \tilde{G}$, which is defined as the kernel of the similitude character on $\tilde{G}$. (Such groups in principle correspond to the group $I$.) We note that there is no moduli problem associated to these groups, hence we cannot follow the path taken by Rapoport and Zink. Instead we rely on methods of Varshavsky \cite{Varshavsky01}, \cite{Varshavsky02}, who proved cases of $p$-adic uniformizations for unitary Shimura varieties by studying their complex uniformization using differential geometry. We first show existence of a canonical model $\SGU$ for the locally symmetric spaces over a field in which $p$ is unramified. This is necessary for the existence of $p$-adic uniformization. We then prove a $p$-adic uniformization result for the varieties $\SGU$ in the spirit of Varshavsky's work. We show that the differential techniques present in his work can be suitably modified to work in our case, which allows us to ultimately prove a level-raising result in which $p$ is not assumed to remain inert in the field $F$.

Then we turn our attention to studying integral structures appearing on spaces of algebraic automorphic forms on these unitary groups. Let $\ell$ be the prime modulo which we seek a level-raising congruence at a place $\nu$ of $F$ split in $E$ such that $I(F_{\nu}) \cong GL_n(F_{\nu})$. We define certain automorphic local systems $V_{\bm{\mu}, \OS}$, $V_{\bm{\mu}, k}$ on the space $\SGU$ where $\OS$ is the ring of integers in an $\ell$-adic field and $k$ is its residue field. The link between the cohomology of these local systems and the spaces of algebraic automorphic forms with $k$-coefficients, $\mathcal{A}(U, W_{\bm{\mu}, k})$, is given by the weight spectral sequence, and is explored by Thorne in his paper \cite{Thorne} for the similitude groups. 

As a second key idea, we first establish a similar relationship for the case of unitary groups.  We use the observation that proving level-raising results is possible if one assumes that the cohomology groups of the $\OS_F[1/\nu]$-arithmetic subgroups of $I(F_{\nu})$ contain no $\ell$-torsion. Thus, we focus on computing the Weight spectral sequence to obtain such torsion-free results. The $E_1$ page of this spectral sequence can be written in terms of algebraic automorphic forms $\mathcal{A}(U, W_{\bm{\mu}, k})$ on the definite unitary group $I$, using their definition as functions on connected components of the special fiber. This is possible using a key combinatorial argument. The $E_2$ page of the weight spectral sequence is related to the cohomology groups of the $\OS_F[1/\nu]$-arithmetic subgroups of $I(F_{\nu})$. Furthermore, this spectral sequence degenerates when $\ell$ is a banal characteristic for $GL_n(E_{\omega})$. (See \cite{Vigneras} for information on banal characteristic.) This follows in particular from the Galois equivariance of the differentials in the spectral sequence. This allows us to prove, combined with $\ell$-torsion-vanishing results of \cite[Theorem 8.12]{LS} (or \cite{Shin} with some further local hypotheses), that the cohomology of the arithmetic subgroups is $\ell$-torsion-free, which in turn suffices to prove level-raising results for algebraic automorphic forms on $I$, as mentioned earlier. This in turn allows us to transfer to $GL_n$ and prove a level-raising result for $GL_n(\A_{E})$. It is worth emphasizing that the hypothesis that $\ell$ is a banal characteristic for $GL_n(E_{\omega})$ plays a crucial role in the argument. 

We now describe the structure of this paper. Section $2$ contains some basic fundamental results concerning the reductive groups occurring in this work. We define the unitary similitude group $\tilde{G}$ over $\Q$ and its unitary subgroup $G$ and their inner forms $\tilde{I}$ and $I$, and set notations. In section $3$, we define the locally symmetric spaces associated to $G$ and use descent from the Shimura varieties for $\tilde{G}$ to describe a model for the locally symmetric spaces over a number field that is unramified at $p$. We then describe $p$-adic uniformizations of these models by using and modifying ideas from \cite{Varshavsky01} \cite{Varshavsky02}. This is the key ingredient that will be used in section $7$. Sections $4$ and $5$ are devoted to recalling the theory of automorphic forms and automorphic systems on the spaces. In section $4$, we define some $\ell$-adic automorphic local systems on these spaces. Then, in section $5$, we define spaces of automorphic forms $\mathcal{A}(U, M)$ and describe the Hecke algebra action upon these spaces, and the cohomology of local systems. Section $6$ is largely for recollection. We recall the definition of the weight spectral sequence associated to a semistable scheme. Finally, section $7$ is where we bring all the arguments together. In this section, we describe the degeneration of the weight spectral sequence for our spaces in characteristic $\ell$ under the banal hypothesis to prove a level-raising result for the group $I$. This involves using a trick involving the Galois action by weights on the terms in the spectral sequence. In the end, in section $8$, we apply this result to deduce level-raising for $GL_n$.

\section{Setup}
\label{sec:Setup}

We set notations and make definitions of the main objects that will be used throughout this paper.

\subsection{Notations and hypotheses}
\label{sec:notation}

We fix an algebraic closure $\overline{\Q}$ of $\Q$ and $n \in \Z_{> 1}$. By $\A_{N}$ we denote the adeles over a number field $N$, and for a finite set of places $S$ of $N$, by $\A_{N}^{S}$ we mean adeles excluding the primes in $S$. $\A_{N, f}$ will denote the finite adeles. If the number field is $\Q$, we drop the subscript $\Q$ from this notation.

Let $E$ be a CM imaginary field with totally real subfield $F$, $[F:\Q] = d$. We assume $d$ is even. Let $\tau'_1, \tau'_2, \ldots, \tau'_d$ be the embeddings at $\infty$ of $F$. Let $D$ be a central division algebra of dimension $n^2$ over $E$. Let $\ast$ be a positive involution of the second kind on $D$ so that the invariants of $\ast$ on $E$ is the subfield $F$. Let $V$ be a left $D$-module and $\psi$ an alternating pairing $\psi : V \times V \rightarrow \Q$ satisfying \begin{equation} \psi(dv, w) = \psi(v, d^{\ast} w) \end{equation} for all $d \in D, v, w \in V$. In fact, we will take $V = D$ and $\psi (x, y) = xy^{\ast}$ from now on.

We define an algebraic group $\tilde{G}$ over $\Q$ by its functor of points : \begin{equation} \tilde{G}(R) := \{ g \in (D \otimes R)^{\times} : g.g^{\ast} \in R^{\ast} \}. \end{equation} where $R$ is a $\Q$-algebra. In particular, $c(g) := g.g^{\ast}$ is the similitude character of $\tilde{G}$.

$\tilde{G}_{\R}$ can be embedded into a product of unitary similitude groups. Namely, for each $i = 1, 2, \ldots, d$, fix an embedding $\tau_i$ of $E$ into $\C$ extending $\tau'_i$. Denote $\Phi:= \{ \tau_1, \ldots, \tau_d \}$.  We can choose an isomorphism $$D \otimes_{\Q} \R \cong \prod_{\tau \in \Phi} D \otimes_{E, \tau} \C \cong \prod_{\tau \in \Phi} M_n(\C)$$ such that $\ast$ corresponds to the operation $ A \rightarrow \overline{A^{t}}$, where $z \rightarrow \overline{z}$ denotes the complex conjugation and $M^{t}$ denotes the transpose of a matrix $M$.

This decomposition induces an orthogonal decomposition with respect to $\psi$ on $$V \otimes_{\Q} \R \cong \prod_{\tau \in \Phi} V \otimes_{E, \tau} \C.$$ There are isomorphisms $ V \otimes_{E, \tau} \C \cong \C^{n} \otimes_{\C} W_{\tau}$ where $M_n(\C)$ acts on the first factor. Then we can find a skew-hermitian form $h_{\tau}$ on $W_{\tau}$ so that \begin{equation} \psi_{\tau}(z_1 \otimes w_1. z_2 \otimes w_2) = tr_{\C / \R}(\overline{z_1}^{t}z_2h_{\tau}(w_1, w_2)). \end{equation} We can find a basis of $W_{\tau}$ such that $h_\tau$ is represented by the matrix \begin{equation} diag(-i, -i, \cdots , -i, i, i, \cdots, i). \end{equation} We denote the number of times $-i$ appears in this matrix by $r_{\tau}$ and the number of times $i$ appears by $r_{\tau^c}$. It can be seen that $r_{\tau} + r_{\tau^c} = (1/n) \dim_{E}V = n$. Thus, the choices we have made so far imply that \begin{equation} \tilde{G}_{\R} \hookrightarrow \prod_{\tau \in \Phi} GU(r_{\tau}, r_{\tau^c}) \end{equation} such that $\tilde{G}_{\R}$ is a normal subgroup with torus cokernel. We now assume that for a chosen embedding $\tau_1$, $r_{\tau_1} = 1$ and for all $\tau \in \Phi$ such that $\tau \not = \tau_1$, $r_{\tau} = 0$.

Under the above identifications, we write $h : \text{Res}_{\C / \R} \G_m \rightarrow \tilde{G}_{\R}$ for the homomorphism \begin{equation} h : z \in \C^{\times} \rightarrow (diag(z, z,\cdots, z, \bar{z}, \bar{z}, \cdots, \bar{z}))_{\tau \in \Phi} \end{equation} where again the number of times $z$ appears is $r_{\tau}$ and the number of times $\bar{z}$ appears is $r_{\tau^{c}}$. Let $X$ denote the $\tilde{G}(\R)$-conjugacy class of $h$ inside the set of homomorphisms $\text{Res}_{\C / \R} \G_m \rightarrow \tilde{G}_{\R}$.

With these assumptions, we have the following result, primarily due to Kottwitz \cite[Page 654]{Kottwitz02} following the theory of canonical models:

\begin{prop}
The pair $(\tilde{G}, X)$ is a Shimura datum. For $\tilde{U} \subset \tilde{G}(\A^{\infty})$ a neat open compact subgroup, the Shimura varieties $S(\tilde{G}, \tilde{U})$ with $S(\tilde{G}, \tilde{U})(\C) : = \tilde{G}(\Q) \backslash \tilde{G}(\A_{f}) \times X / \tilde{U}$ are smooth projective algebraic varieties over $\C$ and admit canonical models over the reflex field $\tau_{1}(Q) \subset \C$.
\end{prop}

Here, the reflex field $\tau_{1}(Q)$ can be explicitly described as the compositum of $\tau_1(E)$ and the reflex field $R$ of the CM-type $\Phi$. (This is shown in section $5.2$, page $134$ of \cite{Clozel}, for instance.)

We define another algebraic group $G$ over $\Q$ by its functor of points : \begin{equation} G(R) := \{ g \in (D \otimes R)^{\times} : g.g^{\ast} = 1 \}. \end{equation} where $R$ is a $\Q$-algebra. There is an exact sequence of algebraic groups \begin{equation}\label{eq:exseq1} 1 \rightarrow G \rightarrow \tilde{G} \xrightarrow{c} \G_m \rightarrow 1. \end{equation}

We let $G'$ denote the derived subgroup of $G$. Note that it is also the derived subgroup of $\tilde{G}$. It is the subgroup of $G$ consisting of matrices of determinant $1$. It is a form of $SL_n$ and in particular, simply connected. There is an exact sequence of algebraic groups \begin{equation} 1 \rightarrow G' \rightarrow G \rightarrow (E^{\times})^{N=1} \rightarrow 1 \end{equation} where $(E^{\times})^{N=1}$ is the torus of elements of $E^{\times}$ with norm $1$ in $F$. We let $\tilde{T} := \tilde{G} / G'$ and $T' := G / G'$. Note that $\tilde{T}$ and $T$ are the largest commutative quotients of the respective groups. We denote the natural map $\tilde{G} \rightarrow \tilde{T}$ by $\nu$. Let $\tilde{Z}, Z$ denote the center of $\tilde{G}, G$ respectively.

Let $p$ be a rational prime. We assume $p$ is unramified in $E$. Let $\nu$ be a place of the reflex field $Q$ whose restriction to $\Q$ is $p$. Choose an isomorphism $\C \cong \C_p$ extending the natural morphism $Q \rightarrow Q_{\nu} \rightarrow \C_p$. The embedding $\tau_1 : E \rightarrow \C$ then corresponds to an embedding $\alpha_1 : E \rightarrow \C_p$ which extends to a continuous embedding $\alpha_1 : E_{\omega_1} \rightarrow \C_p$ for some place $\omega_1$ of $E$ above $p$. Since the group $Gal(Q/E)$ preserves $\tau_1$, the place $\omega$ does not depend on the isomorphism $\C \cong \C_p$. Let $\nu_1$ be the restriction of $\omega_1$ to $F$. To the assumptions made on group $\tilde{G}$, we add the following :

\begin{enumerate}
	\item $\nu_1$ is split in $E$ as $\nu_1 = \omega_{1} \omega_{1}^{c}$.

	\item The Brauer invariant of $D \otimes_{E} E_{\omega_1}$ is $1/n$.
\end{enumerate}

At every other finite place of $E$, $D$ is assumed to be split. (The existence of $\tilde{G}$ satisfying these assumptions follows from \cite[Proposition 2.3]{Clozel} under our assumption that $d$ is even.)

 Let $ \{ \nu_1, \nu_2, \ldots, \nu_t \}$ be the places of $F$ above $p$. We also fix a central simple algebra $\tilde{D}_{\omega_1}$ over $E_{\omega_1}$ of Brauer invariant $1/n$. 

The decomposition $D \otimes_{\Q} \Q_{p} \cong \oplus_{\nu | p} D_{\nu}$ will give an explicit formula for $G(\Q_p)$, namely, $G(\Q_p) \cong \prod_{i=1}^{t} G(\Q_p)_{\nu_i}$ where $G(\Q_p)_{\nu_i}$ is given by $$G(\Q_p)_{\nu_i} = \{ GL_{D_{\nu_i}}(D \otimes_{F} F_{\nu_i}) | \psi(gv, gw) = \psi(v, w) \}. $$ If $\nu_i$ is split in $E$ as $\omega_{i} \omega_{i}^{c}$, we get $G(\Q_p)_{\nu_i} \cong D_{\omega_{i}}^{\times}$ by projection onto first coordinate, and if $\nu_j$ is inert in $E$ giving a place $\omega_j$ of $E$, we get $G(\Q_p)_{\nu_j}$ to be a unitary group at the place $\omega_j$. Consequently, $ G(\Q_p) \cong \tilde{D}_{\omega_1}^{\times} \times G(\Q_p)^{\omega_1} $ where $G(\Q_p)^{\omega_1}$ is a product of general linear groups and unitary groups at $\omega_j \not = \omega_1$ depending upon the splitting of $\nu_j$. In particular, any level subgroup $U_p$ at $p$ splits as $U_{p, \omega_{1}}U_{p}^{\omega_{1}}$.

We now define inner forms of groups $\tilde{G}$ and $G$ obtained by `switching primes'. For this purpose, let $D'$ be a central division algebra of dimension $d^2$ over $E$ equipped with an involution $\ast \ast$ of second kind. We require that $(D, \ast)$ and $(D', \ast \ast)$ are locally isomorphic at all finite places of $E$ except at places above $\nu_1$, where $D'$ is assumed to be split. In addition, we assume $\ast \ast$ to be positive definite at all archimedean places. %We consider $D'$ as a left module over itself and a pairing $\psi' : D' \times D' \rightarrow \Q$ as $\psi'(x, y) = xy^{\ast \ast}$. 
Then we define a reductive group $\tilde{I}$ over $\Q$ by \begin{equation} \tilde{I}(R) := \{ g \in (D' \otimes R)^{\times} : g.g^{\ast \ast} \in R^{\ast} \}. \end{equation} This is an inner form of $\tilde{G}$. We define its unitary subgroup $I$ by \begin{equation} I(R) := \{ g \in (D' \otimes R)^{\times} : g.g^{\ast \ast} = 1 \}. \end{equation} This is an inner form of the group $G$ that is compact at all infinite places. 

For $L$ a finite extension of $\Q_p$, we write $\Omega_L$ for the $(n-1)$-dimensional Drinfeld $p$-adic upper half space over $L$. This is an analytic space obtained by removing from $\mathbb{P}^{n-1}$ the union of all $L$-rational hyperplanes. It is a generic fiber of an explicitly constructed $p$-adic formal scheme $\hat{\Omega}_{L}$, formally locally of finite type over Spf $\OS_L$. There is a faithful action of $PGL_n(L)$ on $\Omega_{L}$ and $\hat{\Omega}_{L}$. We recall that Drinfeld \cite{Drinfeld02} has constructed a projective system $\{ \Sigma^{m}_{L} \}_{m \in \mathbb{N}}$ of $GL_n(L)$-equivariant finite etale Galois coverings of $\Omega_{L} \hat{\otimes} \breve{L}$. We denote by $\Sigma_{L}$ the pro-analytic space $\varprojlim_{m} \Sigma^{m}_{L}$.  We define $$ \mathcal{M}_{L}^{\text{split}} := \hat{\Omega}_{L} \times GL_n(L) / GL_n(L)^{0} $$ where $GL_n(L)^{0} \subset GL_n(L)$ is the open subgroup of matrices with determinant a $p$-adic unit and the sets on the right hand side are identified with the corresponding constant formal schemes over $\OS_L$. We define $\mathcal{M}_{L} := \mathcal{M}_{L}^{\text{split}} \hat{\otimes}_{\OS_L} \OS_{\breve{L}}$. We suppress the field $L$ from the notations when the base field is clear from context.

\section{Canonical models for unitary Shimura varieties and their $p$-adic uniformization}
\label{sec:front-matter}

\subsection{Canonical models}

In this section, we deduce the existence of a canonical model and uniformizations following Deligne/Kottwitz and Varshavsky, for locally symmetric spaces arising from $G$. The idea is to use the existence of these $p$-adically uniformized models, to compute the terms of the weight spectral sequence associated to them in terms of algebraic automorphic forms. Results from this section will be utilized in section \ref{sec:cohomology}.

We use the following hypotheses in this subsection. 
\begin{hyp}
\label{hyp:mild}

For $\tilde{U} $ as above, we set $U := G(\A_f) \cap \tilde{U}$. We assume that for each $\omega_i$ contributing to $G(\Q_p)$ as in section \ref{sec:notation}, we choose the maximal compact subgroup as the level subgroup at that place. We define \begin{equation} \SGU(\C) := G(\Q) \backslash G(\A_f) \times X / U. \end{equation} This is a complex variety. Assume that $\G_{m}(\Q) \cap c(\tilde{U}) = \{1\}$. (This particular mild hypothesis is inspired by a similar hypothesis in Carayol's work \cite[Section 1.4]{Carayol} where he uses the same kind of tricks to give conditions under which a covering map of Shimura curves has the expected Galois group.)
\end{hyp}

There is a natural map $\phi$ from $\SGU(\C)$ to $\tildeSGU(\C)$ induced from the inclusion $i : G \hookrightarrow \tilde{G}$. We show that it is injective.

\begin{lem}
\label{first}
Let $c$ be the similitude character as before. Then under Hypothesis \ref{hyp:mild}, $$\phi : \SGU(\C) \rightarrow \tildeSGU(\C) $$ sending a double coset representative $\overline{(g, x)}$ to $\overline{(g,x)}$ is injective.
\end{lem}
\begin{proof}
Suppose $(g, x), (h, y) \in G(\A_f) \times X $ are such that their image in $\tildeSGU(\C)$ coincides. This implies that there exist $a \in \tilde{G}(\Q), u \in \tilde{U}$ such that  $agu = h \text{ and } ax = y.$ It suffices to prove that $a \in G(\Q)$ and $u \in U$. From the exact sequence (\ref{eq:exseq1}), we see that $$c(agu) = c(au) = c(h) = 1.$$ Thus, $c(a) = c(u^{-1})$ and hence, $c(a) = c(u) = 1$.
\end{proof}

Since $\G_{m}(\Q) \cap c(\tilde{U})$ is a congruence subgroup of $\Z^{\times}$ by definition, we see that the hypothesis we introduce is not very strong. Indeed, a congruence subgroup of $\Z^{\times}$ can be trivial or the group $\{ 1, -1 \}$. So that if $\tilde{U}$ is sufficiently small, then this intersection should be trivial.

We wish to identify $\SGU(\C)$ with its image under $\phi$. This is the content of our next lemma.

\begin{lem}
Under Hypotheses \ref{hyp:mild}, $\phi(\SGU(\C))$ is a union of connected components of the space $\tildeSGU(\C)$ and $\phi$ is an isomorphism onto its image.
\end{lem}
\begin{proof}
This can be seen by writing each space of double cosets as a finite disjoint union of quotients of $X$ by arithmetic subgroups indexed by the set of double coset representatives. ({\it cf. } Milne \cite[Lemma 5.13]{Milne}.)

Alternatively, from the definition of $\phi$, it can be seen that it is a local homeomorphism, and any bijective local homeomorphism is in fact a homeomorphism.
\end{proof}

This identification of connected components can be further explored. There is a zero-dimensional Shimura variety associated with the data of $\tilde{T}$, which can be identified with the set of connected components of $\tildeSGU$. We describe it briefly. See \cite[Chapter 5]{Milne} for more details. Recall that there is a surjective homomorphism $\tilde{Z} \hookrightarrow \tilde{G} \xrightarrow{\nu} \tilde{T}$. We define \begin{equation*} \tilde{T}(\R)^{\dagger} := \text{Im}(\tilde{Z}(\R) \rightarrow \tilde{T}(\R)), \end{equation*} \begin{equation*} \tilde{T}(\Q)^{\dagger} := \tilde{T}(\Q) \cap \tilde{T}(\R)^{\dagger}. \end{equation*} Then, let $$\tilde{C} : = \tilde{T}(\R)/\tilde{T}(\R)^{\dagger}.$$ $\tilde{C}$ is a finite set. For any compact open set $\tilde{K} \subset \tilde{T}(\A_f)$, we can define a zero-dimensional Shimura variety \begin{equation*} S(\tilde{T}, \tilde{K}) := \tilde{T}(\Q) \backslash \tilde{T}(\A_{f}) \times \tilde{C} / \tilde{K}. \end{equation*} It is known that (\cite[Theorem 5.17]{Milne}) - \begin{equation} \pi_{0}(\tildeSGU) \cong S(\tilde{T}, \nu(\tilde{U})). \end{equation}

Similarly, we define \begin{equation} S(T, K) := T(\Q) \backslash T(\A_{f}) \times C / K. \end{equation} where $T(\Q)^{\dagger}$ and $C$ are defined in an analogous fashion. We show that $S(T, \nu(U))$, like its $(G, U)$-counterpart, can be understood as a subgroup of $S(\tilde{T}, \nu(\tilde{U}))$.

\begin{lem}
For objects defined above satisfying hypothesis \ref{hyp:mild}, the natural map $$\theta : S(T, \nu(U)) \rightarrow S(\tilde{T}, \nu(\tilde{U}))$$ is injective.
\end{lem}
\begin{proof}
Note that since $G'$ is contained in the kernel of the similitude character $c$, $c$ is well-defined at the level of $T$. The proof follows the same path we took for Lemma $1$.
\end{proof}

Since $\tilde{T}(\Q)$ is dense in $\tilde{T}(\R)$, $\tilde{C} \cong \tilde{T}(\Q)/\tilde{T}(\Q)^{\dagger}$. Hence, the space $\tilde{T}(\Q) \backslash \tilde{T}(\A_{f}) \times \tilde{C} / \tilde{K}$ can be identified with the set $\tilde{T}(\Q)^{\dagger} \backslash \tilde{T}(\A_{f}) / \tilde{K}$. There is a natural map \begin{equation} \eta : \tilde{G}(\Q) \backslash \tilde{G}(\A_{f}) \times X / \tilde{U} \rightarrow \tilde{T}(\Q)^{\dagger} \backslash \tilde{T}(\A_{f}) / \nu(\tilde{U}) \end{equation} given by $ (g, x) \rightarrow \nu(g) $. There is a similar map $\eta$ induced by $\nu$ from $\SGU(\C)$ to $\STU$. Since the maps $\phi$ and $\theta$ are injective, we can identify $\SGU(\C)$ with a union of the connected components of $\tildeSGU(\C)$ that map to $\STU$ under $\eta$. We collect this information in the proposition below.

\begin{prop}
\label{prop:compo}
The diagram
\begin{equation}
\begin{array}{ccc}

\SGU(\C)  & \hookrightarrow  & \tildeSGU(\C) \\
\downarrow &  & \downarrow \\
\STU & \hookrightarrow & \tildeSTU \\

\end{array}
\end{equation}
is cartesian for $G$, $T$, $U$ and $\tilde{G}$, $\tilde{T}$, $\tilde{U}$ as above.
\end{prop}

 Using the above observation, we show below the existence of canonical models for $\SGU(\C)$ over a field in which $p$ is unramified.

\begin{thm}
\label{Canonical}
For $U$ satisfying Hypothesis \ref{hyp:mild}, there exists a canonical model $\SGU$ of the complex variety $\SGU(\C)$ over $R$ where $R$ is a finite extension of $Q$ unramified at $p$.
\end{thm}
\begin{proof}
We use the fact that $\tildeSGU$ already has a canonical model over $Q$. Milne \cite[Chapter 13]{Milne} describes the explicit action on $\tildeSTU$ of $Aut(\C/E(\tilde{G}, X))$ where $E(\tilde{G}, X)$ is the reflex field. In this case, $E(\tilde{G}, X) = Q$. Since $\tildeSTU/ \STU$ is a finite abelian group, this action factors through a finite abelian quotient of $Aut(\C/Q)$ which implies the existence of a finite abelian extension $R$ of $Q$ over which there exists a model for $\SGU(\C)$. In particular, since the level subgroup $U_p$ is chosen to be the maximal compact subgroup at every place dividing $p$, we conclude that $p$ is unramified in $R$. This is useful for the uniformization result Corollary \ref{main} later, since the Shimura varieties we consider have maximal level at $p$. As such they can be $p$-adically uniformized meaningfully only in this setting by the Drinfeld upper half space.
\end{proof}

In what follows, we denote by $\SGU$ this canonical model.

\subsection{$p$-adic uniformization}

In this section, we describe $p$-adic uniformization of the canonical model $\SGU$ defined in Theorem \ref{Canonical}.  Following Varshavsky \cite{Varshavsky01}, \cite{Varshavsky02}, we set some notation. Since Varshavsky works with restrictions of groups defined over $F$ instead of directly working with groups defined over $\Q$, the notation in this section is slightly different. We explain the connection between the two versions after setting notation here.

In this section, all level subgroups that appear are assumed to be sufficiently small in the usual sense, for algebraization and other such useful notions to hold for our Shimura varieties.

Let $(D, \ast)$ be as before. Then one can define $\Gint : = \mathbf{GU}(D, \ast)$, an algebraic group over $F$ of unitary similitudes, as a functor of points -

$$ \Gint(R) = \{ d \in (D \otimes_{F} R)^{\times} | d . d^{\ast} \in R^{\times} \} $$ for every $F$-algebra $R$. Then we define $\Hint = Res_{F / \Q} \Gint$, where $Res_{F/ \Q}$ denotes the Weil restriction from $F$ to $\Q$. The homomorphism $h : Res_{\C / \R} \mathbb{G}_m \rightarrow \Hint(\R)$ is defined in \cite{Varshavsky02}, in the case of our signature, by letting for each $z \in \C^{\times}$, $$ h(z) = (diag(1, \cdots, 1, z/\bar{z})^{-1}; I_n; \cdots; I_n).$$ Let $\mathbf{M^{int}}$ be the conjugacy class of $h$. %It follows by comparing this definition with the definition of our conjugacy class, that our uniformizations will differ by the corresponding twist, and that it suffices to prove a uniformization result for either of the conjugacy classes. (See remark $3.13$ of \textit{loc. cit.}) 

Let $\mathbf{\tilde{X}^{int}}$ be the canonical model associated to the Shimura variety for the Shimura datum $(\Hint, \mathbf{M^{int}})$. $\mathbf{\tilde{X}^{int}}$ is an inverse limit over level subgroups $T$ of the quotients $\mathbf{\tilde{X}^{int}} / T$. Then, $\mathbf{\tilde{X}^{int}}$ is a scheme defined over the reflex field of this Shimura datum, say $Q'$. $Q'$ is described explicitly on page $81$ of \cite{Varshavsky02}. In particular, if $v$ is a prime above $p$ of $Q'$, we set $X^{int} := \mathbf{\tilde{X}^{int}} \otimes_{Q'} Q'_{v}$. Then, the First Main Theorem of \textit{loc. cit.} describes a $p$-adic uniformization of $X^{int}$ by the Drinfeld upper half space as a scheme over $Q'_{v}$. We describe this theorem briefly below, since we will use it to deduce our uniformization result.  %Recall that there exists a canonical model for the locally symmetric spaces $\SGU \subset \tildeSGU$ defined over an extension of $F$, say $L$, in which $p$ is unramified. Let $\nu$ be a prime of $L$ lying above $p$ such that its restriction to $F$, $\nu_{1}$, splits in $E$ as $\nu_1 = \omega \omega^{c}$ and $\omega, \omega^{c}$ are the only places where the division algebra $D$ is ramified with invariants $1/n$ and $-1/n$. Then we set $$X^{int} : = \mathbf{\tilde{X}^{int}} \otimes F_{\nu_{1}}, $$ and
\begin{rem}
\label{rem:switch}
Before we begin, we discuss the connection of group $\mathbf{H^{int}}$ with the group $\tilde{G}$. Indeed, Rapoport and Zink \cite{RZ}, for instance, work as we did in the previous section with Shimura data associated to the group $\tilde{G}$ as opposed to the group $\mathbf{H^{int}}$. The former has similitude character valued in $\mathbb{G}_m / \mathbb{Q}$ while the latter has similitude character valued in $Res_{F/ \mathbb{Q}} \mathbb{G}_m$. If we denote the analogous canonical model associated to the Shimura variety arising from $\tilde{G}$ as a scheme $\mathbf{(\tilde{X}^{int})'}$, then $\mathbf{(\tilde{X}^{int})'}$ is an abelian unramified twist of $\mathbf{\tilde{X}^{int}}$, possibly defined over a different reflex field $Q$ (as denoted earlier). This difference between $Q$ and $Q'$ can be explicitly computed by computing the Galois action on the set of special points. In particular, it is proven in Lemma $2.6$ of \cite{Varshavsky02} that $Q'_{v} \cong E_{\omega_1}$. Further, the difference between Galois actions on $\mathbf{(\tilde{X}^{int})'}$ and $\mathbf{(\tilde{X}^{int})}$ is given by some explicit homomorphism $Aut(\C / Q') \rightarrow F_{\nu_1}^{\times}$ as shown in Remark $3.13(a)$ of \textit{loc. cit}.  

So, as mentioned in Remark $3.13(b)$ of \textit{loc. cit.}, it is sufficient to prove that the scheme $X^{int}$ admits a $p$-adic uniformization by the Drinfeld upper half space over $Q'_{v}$ to deduce by twisting that the corresponding scheme for $\mathbf{(\tilde{X}^{int})'}$ admits a $p$-adic uniformization over $Q_{v'}$ for $v'$ a corresponding place over $p$ of $Q$. Finally, since we work with unitary subgroup $G$ of $\tilde{G}$, proving a $p$-adic uniformization result for the corresponding unitary subgroup of $\mathbf{H^{int}}$ (over its field of definition) is sufficient to deduce the $p$-adic uniformization for the space arising from $G$ (over its corresponding field of definition) by a similar argument.
\end{rem}
 
For the purpose of describing this $p$-adic uniformization, we first have to describe an inner form of $\Hint$ obtained by switching primes. This is given using the central division algebra $D'$ over $E$ and an involution of second kind $(D', \ast\ast)$ defined earlier. Then we define $\mathbf{\tilde{I}} : = \mathbf{GU}(D', \ast \ast)$ as a group over $F$ analogously to the definition of $\mathbf{G^{int}}$. This is an inner form of $\mathbf{G^{int}}$. Recall that we have also fixed a central simple algebra $\tilde{D}_{\omega_1}$ over $E_{\omega_1}$ of invariant $1/n$. %Then we define $$ \mathcal{I'} : = F_{\nu_1}^{\times} \times \mathbf{I}(\A_{F}^{f;\nu_1}) $$ % 

%and $$ \tilde{\mathcal{G}} : = \tilde{D_{\omega}}^{\times} \times \mathcal{I'}.$$% 

Then we have isomorphisms $$\mathbf{\tilde{I}}(\A_{F, f}) \cong \GLn(E_{\omega_1}) \times F_{\nu_1}^{\times} \times \mathbf{\tilde{I}}(\A_{F, f}^{\nu_1}) $$ and $$ \Gint(\A_{F, f}) \cong  \tilde{D}_{\omega_1}^{\times} \times F_{\nu_1}^{\times} \times \mathbf{\tilde{I}}(\A_{F, f}^{\nu_1}). $$ 

For each compact open subgroup $ V \subset F_{\nu_1}^{\times} \times \mathbf{\tilde{I}}(\A_{F, f}^{\nu_1})$, one considers the double quotient $$ X_{V, m}^{an} : = \mathbf{\tilde{I}}(F) \backslash \Sigma_{E_{\omega_1}}^{m} \times F_{\nu_1}^{\times} \times \mathbf{\tilde{I}}(\A_{F, f}^{\nu_1}) / V $$ where $\Sigma_{E_{\omega_1}}^{m}$ is the finite etale Galois covering of the rigid analytic Drinfeld $p$-adic upper half space over $E_{\omega_1}$ from section \ref{sec:notation}. For every $V$ as above, Proposition $2.8$ of \cite{Varshavsky02} shows that this has a structure of an analytic space over $E_{\omega_1}$, and in fact it arises from a projective scheme $X_{V, m}$ over $E_{\omega_1}$. Set $$X : = \varprojlim_{V, m} X_{V, m}.$$ Then, the group $\mathbf{G^{int}}(\mathbb{A}_{F, f})$ acts on $X$ in such a way, that for each level subgroup $T$ of $\mathbf{G^{int}}(\mathbb{A}_{F, f})$, $X / T$ is a projective scheme, and $X \cong \varprojlim_{T} X / T$. (See \cite[Con. 1.5.1]{Varshavsky01} for a justification.) Then following is the First Main Theorem of \cite{Varshavsky01}, \cite{Varshavsky02} :

\begin{thm}[Varshavsky]
\label{Varsha}
There exists a $\Gint(\A_{F, f})$-equivariant isomorphism of schemes over $E_{\omega_1}$ $$ \varphi : X \xrightarrow{\sim} X^{int}. $$
\end{thm}

%\begin{rem}
%\begin{enumerate}
%\item Note that $F_{\nu_1}$ and $E_{\omega}$ are in fact the same fields. Hence both objects are defined over the same base.
%
%\item Unraveling the definitions of above groups via Weil restriction, we see that proving $p$-adic uniformization results for spaces above (a priori defined over $F$) gives $p$-adic uniformization results for spaces over $\Q$.
%\end{enumerate}
%\end{rem}

We will use this result to prove a $p$-adic uniformization result for unitary Shimura varieties arising from a unitary subgroup of $\mathbf{G^{int}}$. (These, as noted in Remark $1$, are a twist of $S(G, U)$.)

\begin{rem}
\label{rem:level}
We briefly explain the appearance of $\Sigma_{E_{\omega_1}}^{m}$. Let $T_r$ denote the $r$-th congruence subgroup of $\OS_{\tilde{D}_{\omega_1}}$. Then, the Shimura variety corresponding to the level subgroup $T_r \times V \subset \mathbf{G^{int}}(\A_{F}^{f})$ is uniformized by $\Sigma^{r}_{E_{\omega_1}}$ as explained in the introduction (page $58$) and Main Theorem of \cite{Varshavsky01}. It is technically better to work with the inverse limit over all $r$, hence the definition and the statement of theorem as stated. In our particular case, we have $r = 1$ and hence we obtain uniformization by $\Sigma^{1}_{E_{\omega_1}} \cong \Omega_{E_{\omega_1}}$.
\end{rem}

Let $\mathbf{U^{int}}$ be the unitary subgroup of $\mathbf{G^{int}}$. Then, the Shimura variety associated to $Res_{F/ \Q}\mathbf{U^{int}}$ has a canonical model $\mathbf{Y^{int}}$ which is a scheme over some reflex field, say $R'$, which we take to be large enough to contain $Q'$. $\mathbf{Y^{int}}$ is also an inverse limit over all level subgroups. Let $u$ be a place of $R'$ over $v$ and let $Y^{int} := \mathbf{Y^{int}} \otimes_{R'} R'_{u}$. 

We define $\mathbf{I}$ to be the unitary subgroup of $\mathbf{\tilde{I}} = \mathbf{GU}(D, \ast\ast)$. Then, as before, for compact open subgroup $W \subset \mathbf{I}(\mathbb{A}^{\nu_1}_{F, f})$, we define the quotients $$ Y_{W, m}^{an} :=  \mathbf{I}(F) \backslash (\Sigma_{E_{\omega_1}}^{m} \otimes_{E_{\omega_1}} R'_{u}) \times \mathbf{I}(\mathbb{A}_{F, f}^{\nu_1}) / W. $$ Then, for every W, this has the structure of an analytic space, and arises from a projective scheme $Y_{W, m}$ over $R'_{u}$. We set $$Y := \varprojlim_{W, m} Y_{W, m}.$$ Then, for each level subgroup $T'$, $Y / T'$ has the structure of a projective scheme and $Y \cong \varprojlim_{T'} Y / T'$.

%Then we set $$I' := \mathbf{I}(\A_{F}^{f; \nu_1})$$ and $$\tilde{I} : = \tilde{D}_{\nu_1}^{\times} \times I'$$ where $\tilde{D}$ is as before. There is a construction of an analytic space which is projective, say $Y$ for this true unitary group, completely analogous to the construction of $X$ in the unitary similitude group case, viz. it is given, for a level subgroup $U'$, by : $$ (U' \backslash Y)^{an} \cong U' \backslash \Omega_{F_{\nu_1}} \times I' / \mathbf{I}(F). $$ and $Y : = \varprojlim Y_{U'}.$

We now come to the main theorem in this section.

\begin{thm}
\label{Uniformization}
For $Y, Y^{int}$ schemes over $R'_{u}$ as above, there exists an $\mathbf{U^{int}}(\mathbb{A}_{F, f})$-equivariant isomorphism $$ \varphi_{I, G} : Y \xrightarrow{\sim} Y^{int}.$$
\end{thm}

 One can obtain a proof by suitably modifying techniques in \cite{Varshavsky02}, by proving that, in the notation that appears there, the stabilizer $\Delta$ of a connected component $M$ is isomorphic to $PU(n-1, 1)$. We however prefer to give an alternate proof inspired by the example of quaternionic Shimura varieties.

\begin{proof}
We denote by $\tilde{\mathcal{G}} := \tilde{D}_{\omega_1}^{\times} \times F_{\nu_1}^{\times} \times \mathbf{\tilde{I}}(\A_{F, f}^{\nu_1}) = \mathbf{G^{int}}(\mathbb{A}_{F, f})$. We also denote by $\tilde{\mathcal{U}} := \tilde{D}_{\omega_1} \times \mathbf{I}(\mathbb{A}_{F, f}^{\nu_1}) = \mathbf{U^{int}}(\mathbb{A}_{F, f})$. We have the embeddings of algebraic groups $ \mathbf{I} \hookrightarrow \mathbf{\tilde{I}} $. This induces embedding of groups $$ \tilde{\mathcal{U}} = \tilde{D}_{\omega_1}^{\times} \times \mathbf{I}(\A_{F, f}^{\nu_1}) \hookrightarrow  \tilde{D}_{\omega_1}^{\times} \times  \mathbf{\tilde{I}}(\A_{F, f}^{\nu_1}) = \mathcal{G}_1 \backslash \mathcal{\tilde{G}}$$ where $\mathcal{G}_1$ is simply the group $ \{ 1 \} \times F_{\nu_1}^{\times} \times \{ 1 \} \subset \tilde{\mathcal{G}}$, which keeps track of similitude factors. (See \cite[Lemma 5.4]{Varshavsky02} for an explanation of the appearance of this nonstandard quotient.)

Similarly,  the embedding $ \mathbf{U^{int}} \hookrightarrow \Gint $ of algebraic groups induces an embedding $$ \mathbf{U^{int}}(\A_{F, f}) \hookrightarrow \Gint(\A_{F, f}) \rightarrow \mathcal{G}_1 \backslash \Gint(\A_{F, f}) $$ commuting with the above identifications.

Then, we see that the natural equivariant embeddings $$ \mathbf{I}(F) \backslash ( \Sigma_{E_{\omega_1}} \times \mathbf{I}(\A_{F, f}^{\nu_1}) )  \hookrightarrow \mathbf{\tilde{I}}(F) \backslash (\Sigma_{E_{\omega_1}} \times \mathbf{\tilde{I}}(\A_{F, f}^{\nu_1}) )  $$ and $$ \mathbf{U^{int}}(F) \backslash  (\mathbf{M^{int}} \times \mathbf{U^{int}}(\A_{F, f}) ) \hookrightarrow \Gint(F) \backslash ( \mathbf{M^{int}} \times (\mathcal{G}_1 \backslash \Gint(\A_{F, f}) ) $$ induce, by GAGA, an embedding $Y \hookrightarrow \mathcal{G}_1 \backslash X$ and an embedding $ i : Y^{int}_{\C} \hookrightarrow \mathcal{G}_1 \backslash X_{\C}^{int}$. We then note that $i$ is $R'_{u}$-rational, by a computation essentially the same as our earlier computation of the field of definition of $\SGU$. A similar computation is made in great detail in \cite[Step 4, Section 3]{Varshavsky02}. Thus, we get an embedding of schemes $i : Y^{int} \hookrightarrow \mathcal{G}_1 \backslash X^{int}$. We treat these embeddings as subsets.

We note now, that $  \mathbf{PI}  \xrightarrow{\sim} \mathbf{P\tilde{I}} $ and $ \mathbf{PU^{int}} \xrightarrow{\sim} \mathbf{PG^{int}} $. This implies that modulo centers, above embeddings of schemes are actually isomorphisms : $$ Z(\tilde{\mathcal{U}}) \backslash Y \cong Z(\tilde{\mathcal{G}}) \backslash X $$ and $$Z(\mathbf{U^{int}}(\A_{F, f})) \backslash Y^{int} \cong Z(\Gint(\A_{F, f})) \backslash X^{int}.$$

By Theorem \ref{Varsha}, we have an isomorphism equivariant for $\tilde{\mathcal{G}}$-action $$\varphi : X \xrightarrow{\sim} X^{int}.$$ It induces an isomorphism, equivariant for $\mathcal{G}_1 \backslash \tilde{\mathcal{G}}$-action $$ \tilde{\varphi} : = \mathcal{G}_1 \backslash X \xrightarrow{\sim} \mathcal{G}_1 \backslash X^{int}. $$

Let $x \in Y \subset \mathcal{G}_1 \backslash X$. Then there exists $g \in Z(\tilde{\mathcal{G}})$ such that $$ g \tilde{\varphi}(x) \in Y^{int} \subset \mathcal{G}_1 \backslash X^{int}.$$ Thus, $g \tilde{\varphi}$ maps the $\tilde{\mathcal{U}}$-orbit of $x$ into $Y^{int}$. By \cite[Proposition 1.3.8 and 1.5.3]{Varshavsky01}, the $\mathbf{U^{int}}(\A_{F, f})$-orbit is Zariski dense in both $Y$ and $Y^{int}$, we see that $g\tilde{\phi}$ restricts to the required isomorphism $ \varphi_{I, G} : Y \xrightarrow{\sim} Y^{int}.$

\end{proof}

We record the following corollary, which is more suited for application to the problem of level-raising. This is a statement analogous to a consequence of corollary $6.51$ of \cite{RZ}, which is a $p$-adic uniformization result of Rapoport-Zink. We now drop the subscript from the notation of the place $\nu_1$ and refer to it as simply a fixed place $\nu$. Similarly we denote $\omega_1$ by $\omega$. We write the subgroup $U_{p, \omega}$ of $U$ as $U_{\omega}$ now, and denote by $U^{\omega}$ the subgroup $U^{p}U^{\omega}_{p}$ and write $U = U^{\omega}U_{\omega}$. We denote by $G(\A^{\omega, \infty})$ the open compact subgroup of $G(\A_f)$ defined in the obvious way by excluding contribution of the $\omega$ factor at $p$. Recall from section \ref{sec:notation} that we have : 

\begin{equation} \mathcal{M}^{\text{split}}_{E_{\omega}} := \hat{\Omega}_{E_{\omega}} \times GL_n(E_{\omega}) / GL_n(E_{\omega})^{0}. \end{equation}  Analogously we have $\mathcal{M}_{E_{\omega}} := \mathcal{M}^{\text{split}}_{E_{\omega}} \hat{\otimes}_{\OS_{E_{\omega}}}\OS_{\breve{E_{\omega}}}$. Then, we have, as a corollary -

\begin{cor}
\label{main}
For canonical models $\SGU$ satisfying Hypothesis \ref{hyp:mild}, there exists a model $Sh(G, U)$ over $\OS_{R_{u}}$ of $S(G, U)$ and an isomorphism -

$$ I(\Q) \backslash [\mathcal{M}_{E_{\omega}} \times G(\A^{\omega, \infty}) / U^{\omega} \cong (Sh(G, U) \otimes_{\OS_{R_{u}}} \OS_{\breve{R_{u}}})^{\wedge} $$ which follows from Theorem \ref{Uniformization}. 
\end{cor}

\begin{proof}
We recall from remark \ref{rem:switch} that the objects associated to $\tilde{G}$ and $G$ are twists of the objects associated to $\Hint$ and its unitary subgroup. So the chief issue is to see how to associate to quotients of $Y^{int}$ an integral structure from the $p$-adic uniformization, since then the corresponding twists will acquire it as well, as mentioned in Remark $3.13(b)$ of \cite{Varshavsky02}. We also recall that even as $Y$ contains all Drinfeld \'{e}tale coverings of the Drinfeld upper half plane, we only need focus on $\Sigma^{1} \cong \Omega$ for uniformization of our Shimura varieties, since the level is maximal at $p$. For a reference and a precise statement, look at \cite[Theorem 3.2.1, Corollary 3.2.2]{Varshavsky01}. (In particular, for a level subgroup of the form $ S = T_1 \times V$ as in Remark \ref{rem:level}, one can deduce from $ \varphi_{I, G} : Y \xrightarrow{\sim} Y^{int}$ an isomorphism $(Y^{int}/S)^{an} \cong GL_{d}(R'_{u}) \setminus (\Omega_{R'_{u}} \times (V \backslash \mathbf{I}(\A_{F, f}) / I(F) ) )$.)

Let us base change to $\breve{R}'_{u}$. Then, $\Omega_{\breve{R}'_{u}}$ is an analytic space, and in fact it is the generic fiber of the $p$-adic formal scheme $\hat{\Omega}$ over $\OS_{\breve{R}'_{u}}$, such that the action of $PGL_{n}$ extends to this formal scheme. Thus, there is a natural structure of a formal $\OS_{R_{u}}$-model $\mathcal{Y}_S$, say, arising via quotient of the formal scheme $\hat{\Omega}$, on the quotient $(Y^{int}/S)^{an}$ via the isomorphism induced by $\varphi_{I, G}$ on generic fibers and thus we have the required structure. 

Furthermore, the scheme $\hat{\Omega}$ comes with a natural descent datum described in \cite[Corollary 6.44]{RZ}. We use this natural descent datum on $Y^{int}/S$ to define an integral model $Sh(G, U)$ as required. 

%Indeed theorems \ref{Varsha} and \ref{Uniformization} are a priori isomorphisms of schemes over $Q'_{v}$ and $R'_{u}$ respectively. Since the $p$-adically uniformized spaces $X$ and $Y$ themselves have formal $\OS_{Q'_{v}}$ and $\OS_{R'_{u}}$ models, we can use the isomorphisms in these theorems to define the integral structure to the schems $X^{int}$ and $Y^{int}$ taking the twisting mentioned in Remark $1$ into account. The isomorphism follows. 

(This argument is unravelling of a comment in remark $3.13(b)$ of \cite{Varshavsky02} that \cite[Theorem 6.50, Corollary 6.51]{RZ} follows from the main theorem there after the twisting that we mentioned. See also \cite[Lemma 3.2 ]{Thorne} for a discussion on quotients of $\Omega$.) 
\end{proof}

\begin{rem}
We note that the map in the statement of the corollary no longer has a moduli interpretation, which is used in \cite{RZ}. It is instead the ``union of connected components'' incarnation that we have mentioned earlier. It might be possible to give a moduli interpretation in this case, by putting an appropriate condition on the determinant of the polarization that appears in the classical moduli problem, but we have not pursued such an approach here.
\end{rem}

\section{Automorphic local systems}
\label{sec:Automorphic local systems}

From now on we will only consider sufficiently small open compact subgroups $U $ as in Hypothesis \ref{hyp:mild} so that Corollary \ref{main} holds true. We now describe some local systems on $\SGU$ that are related to algebraic representations of $G$ and spaces of automorphic forms.

Corresponding to the chosen CM type $\Phi$, we get an isomorphism $$ G(\C) \cong \prod_{\tau \in \Phi} GL_n(\C).$$ Let $T \subset G \otimes_{\Q} \C$ be the product of diagonal maximal tori : $$ T(\C) \cong \prod_{\tau \in \Phi} \C^{\times} \times \ldots \times \C^{\times}.$$ Then we can write $ X^{\ast}(T) \cong (\Z^{n})^{\Phi}$ and we denote by $X^{\ast}(T)_{+}$ the subset of dominant weights $\bm{\mu} = (\mu_{\tau})_{\tau \in \Phi}$, which are precisely the ones satisfying the condition $$ \mu_{\tau, 1} \ge \mu_{\tau, 2} \ge \ldots \ge \mu_{\tau, n} $$ for each embedding $\tau$ in $\Phi$. If $\ell$ is a rational prime, we say $\bm{\mu}$ is $\ell$-small if, for each $\tau \in \Phi$, $$ 0 \le \mu_{\tau, i} - \mu_{\tau, j} < \ell $$ for all $i, j$ such that $0 \le i < j \le n$. For $\ell$ unramified in $E$ and $\bm{\mu}$ $\ell$-small, we can associate to $\bm{\mu}$ an $\ell$-adic local system on $\SGU$ as follows. (See \cite{HT}, section $III.2$, \cite{Harris}, section $7.1$.)

We fix an isomorphism $ \iota : \bar{\Q}_{\ell} \cong \C $. Let $K$ be a finite extension of $\Q_{\ell}$ in $\bar{\Q}_{\ell}$ such that by assumption, the algebraic representation of $G \otimes_{\Q} \Q_{\ell}$ of highest weight $\iota^{-1} \bm{\mu}$, say $W_{\bm{\mu}, K}$, is defined over $K$. Let $\mathscr{O}$ be the ring of integers of $K$ with maximal ideal $\lambda$, and residue field $k$. Let $U_{\ell}$ be a hyperspecial maximal compact subgroup of $G(\Q_{\ell})$. Then, up to homothety, there is a unique $U_{\ell}$-invariant $\mathscr{O}$-lattice of $W_{\bm{\mu}, K}$. We choose one such and call it $W_{\bm{\mu}, \mathscr{O}}$. Its uniqueness follows from the $\ell$-small hypothesis by considering the reduction modulo $\ell$.

Then, given an integer $m \ge 1$, we define $U(m) \subset U$ be an open compact subgroup acting trivially on $W_{\bm{\mu}, \mathscr{O}/\lambda^{m}} = W_{\bm{\mu}, \mathscr{O}} \otimes_{\mathscr{O}} \mathscr{O}/\lambda^{m}$ and write it as $U(m) = U^{\omega}(m)U_{\omega}$. Then $U$ acts on the constant sheaf defined by $W_{\bm{\mu}, \mathscr{O}/\lambda^{m}}$ on $S(g, U(m))$, and the quotient defines an \'{e}tale local system on $\SGU$, which we denote by $V_{\bm{\mu}, \mathscr{O}/\lambda^{m}}$. We then take $$V_{\bm{\mu}, \mathscr{O}} : = \varprojlim_{m} V_{\bm{\mu}, \mathscr{O}/\lambda^{m}} $$ and $$V_{\bm{\mu}, K} : =  V_{\bm{\mu}, \mathscr{O}} \otimes_{\mathscr{O}} K.$$ We note that sections of the local system $V_{\bm{\mu}, \mathscr{O}/\lambda^{m}}$ over an \'{e}tale open $T \rightarrow \SGU$ can be identified with the set of functions $f : \pi_{0}(S(G, U(m))\times_{S(G, U)} T) \rightarrow W_{\bm{\mu}}$ such that for all $\sigma \in U, C \in \pi_{0}(S(G, U(m))\times_{S(G, U)} T)$, we have $f(C \sigma) = \sigma^{-1}f(C).$

\section{Spaces of Automorphic forms}
\label{sec:Spaces of Automorphic forms}

We recall the definition of spaces of automorphic forms with integral coefficients on the definite unitary group $I$. We recall the Hecke action on these spaces and the definition of cohomology of these spaces.

We take $\ell$, $K$, $\mathscr{O}$, and $k$ as before. Let $U_{\ell} \subset I(\Q_{\ell})$ be an open compact subgroup, and suppose $M$ is a finite $\mathscr{O}$-module on which $U_{\ell}$ acts continuously. Then the space of automorphic forms $\mathcal{A}(M)$ is defined to be the set of locally constant functions $f : I(\A^{\infty}) \rightarrow M$ such that for all $\gamma \in \Q, f(\gamma g) = f(g)$. This can be endowed with an action of the group $I(\A^{\infty}) \times U_{\ell}$ by the formula $(g.f)(h) := g_{\ell}f(hg)$ where $g_{\ell}$ is the projection to the $\ell$-th component. If $U \subset I(\A^{\infty}) \times U_{\ell}$ is a subgroup, then we define $\mathcal{A}(U, M) := \mathcal{A}(M)^{U}$.

Then we record the following lemma on admissibility of these as representations (Lemma $2.2$ from \cite{Thorne})-

\begin{lem}
For $p \neq \ell$ prime, let $U^{\omega}$ be an open compact subgroup of $I(\A^{\omega, \infty})$ whose projection to $I(\Q_{\ell})$ is contained in $U_{\ell}$. Then for any subgroup $U_{p, \omega} \subset I(\Q_p)$ at the place $\omega$, $\mathcal{A}(U^{\omega}, M)^{U_{p, \omega}}$ is a finite $\mathscr{O}$-module.
\end{lem}

Let $\bm{\mu}$ be a choice of $\ell$-small dominant weight, and let $U \subset I(\A^{\infty})$ be an open compact subgroup. Then there is a finite free $\mathscr{O}$-module $W_{\bm{\mu}, \mathscr{O}}$ and a corresponding space of automorphic forms $\mathcal{A}(U, W_{\bm{\mu}, \mathscr{O}})$, which is a finite free $\mathscr{O}$-module. This space has an interpretation as an isotypic component of the space of automorphic forms $\mathcal{A}$ on $I$ as follows: Let $W_{\bm{\mu}, \C}$ denote the representation of $I(\R) \subset I(\C)$ which is the restriction of the algebraic representation with highest weight $\bm{\mu}$. Then we have an isomorphism: $$ \mathcal{A}(U, W_{\bm{\mu}, \mathscr{O}}) \otimes_{\mathscr{O}, \iota} \C \cong Hom_{I(\R)}(W^{\vee}_{\bm{\mu}, \C}, \mathcal{A}).$$

\subsection{Hecke action}
If $T$ is a finite set of primes containing $\ell$ and not containing $p$, such that $U_q$ is a hyperspecial maximal compact subgroup for all $q \not \in T$, let $\T_{T}^{univ} = \mathscr{O}[ \{ T_{1}^{\nu}, \ldots, T_{n}^{\nu}, (T_{n}^{\nu})^{-1} \} ]$ denote the universal Hecke algebra in infinitely many variables corresponding to the unramified Hecke operators at places $\nu$ of $F$ which split in $E$ and are away from $T$. This Hecke algebra acts on $\mathcal{A}(U, W_{\bm{\mu}, \mathscr{O}})$ by $\mathscr{O}$-algebra endomorphisms, and on the cohomology of the local systems, $H^{i}(S(G, U)_{\overline{E}}, W_{\bm{\mu}, \OS})$ via the isomorphism $I(\A^{p, \infty}) \cong G(\A^{p, \infty})$.

If $\sigma$ is an automorphic representation of $I(\A)$ such that $(\sigma^{\infty})^{U} \neq 0$ and $\sigma_{\infty} \cong W^{\vee}_{\bm{\mu}, \C}$, we can associate to it a maximal ideal $m_{\sigma} \subset \T^{univ}_{T}$ by assigning to each Hecke operator the reduction modulo $\ell$ of its eigenvalue. If $\sigma^{'}$ is another automorphic representation of $I(\A)$, we say that $\sigma^{'}$ {\textit contributes to} $ \mathcal{A}(U, W_{\bm{\mu}, \mathscr{O}})_{m_{\sigma}} $ if  $((\sigma^{'})^{\infty})^{U} \neq 0$ and $\sigma^{'}_{\infty} \cong W^{\vee}_{\bm{\mu}, \C}$, and the intersection of $\iota^{-1}((\sigma^{'})^{\infty})^{U}$ and $ \mathcal{A}(U, W_{\bm{\mu}, \mathscr{O}})_{m_{\sigma}} $ inside $ \mathcal{A}(U, W_{\bm{\mu}, \mathscr{O}}) \otimes_{\OS} \bar{\Q}_{\ell} $ is nontrivial.

\subsection{Cohomology}
For $G = GL_n(L)$ where $L$ is a finite extension of $\Q_p$ and $\OS$ as before, cohomology groups $H^{\ast}(M)$ of a smooth $\OS[G]$-module $M$ are defined as follows.
 
Let $U_0 = GL_n(\OS_{L})$ be the standard maximal compact subgroup and $B \subset U_0$ be the standard Iwahori subgroup. Let $\zeta$ denote the matrix

$$  \left( \begin{array}{ccccc}
0 & 1  & 0 & \dots & 0 \\
0 & 0 & 1 & \dots &  0 \\
\vdots & \vdots &  & \ddots & \vdots \\
0 & 0 & \dots & 0 & 1 \\ 
\varpi & 0 & \dots & 0 & 0 \end{array} \right)  $$

where $\varpi$ is a fixed uniformizer for $\OS_L$. For $ i = 0, \ldots, n-1$, let $U_i = \zeta^{-i}U_{0}\zeta^{i}$. All of these contain $B$. For $A \subset \{ 0, 1, \ldots, n-1 \}$, we write $U_{A} : = \cap_{i \in A} U_{i}$. Then we define a complex $C^{\bullet}(M)$ by the formula $$ C^{i}(M) = \bigoplus_{A \subset \{0, 1,  \cdots, n-1 \}} M^{U_A}, $$ where the direct sum is over subsets of cardinality $i+1$. The differential $d_i : C^{i} \rightarrow C^{i+1}$ is given by the sum of the restriction maps $r_{A, A^{'}} : M^{U_A} \rightarrow M^{U_{A^{'}}}$ for $A \subset A^{'}$, each multiplied by a sign $\epsilon(A, A^{'})$ which is given as follows. If $A^{'} = \{ i_1, \cdots, i_r \}$ written in increasing order, and $A = A^{'} \setminus \{ i_s \}$, then \begin{equation} \label{sign} \epsilon(A, A^{'}) = (-1)^{s}.\end{equation} Then $H^{\ast}(M)$ is defined to be the hypercohomology of this complex.

The cohomology groups $H^{i}(\mathcal{A}(U^{\omega}), W_{\bm{\mu}, k})$ will be key for proving the level-raising result. 

\section{Weight spectral sequence}
\label{sec:Weight spectral sequence}

We recall the weight spectral sequence of Rapoport and Zink. A standard reference for this section is \cite{Saito} and we largely follow it. No originality is claimed. We need to compute this sequence for our Shimura varieties, and we will use the $p$-adic uniformizations of these varieties to this end. For now, we recall basic properties of this sequence. In particular, the fact that this sequence can be understood as the spectral sequence associated to a filtered object is crucial for its application in the later sections.

Let $S : = Spec \OS_{F_0}$ be the spectrum of a complete discrete valuation ring $\OS_{F_0}$. As usual, $s$ is the closed point of $S$ and $\eta$ is the generic point of $S$. Let $F_0 := Frac \OS_{F_0}$, and we fix an algebraic closure $\bar{F_0}$ of $F_0$. We write $\bar{s}, \bar{\eta}$ for the geometric points of $S$ above $s, \eta$ respectively. Suppose that $ f: X \rightarrow S $ is a proper, strictly semistable morphism of relative dimension $n$. Then the special fiber $X_s$ is a strict normal crossings divisor on $X$. We write $X_1, \ldots, X_h$ for its irreducible components. We suppose that each irreducible component $X_i$ is smooth over $\kappa(s)$. For $A \subset \{ 1, \ldots, h \}$ we write $X_A$ for the intersection $\bigcap_{i \in A} X_{i}$, and $ X^{(m)} : = \bigsqcup_{|A| = m+1} X_{A}$ for the disjoint union. Let $K$ be a finite extension of $\Q_{\ell}$ with ring of integers $\OS$, uniformizer $\lambda$, and residue field $k$, where $\ell$ is coprime to the residue characteristic of $\OS_{F_0}$. Let $\Lambda = K, \OS$, or $k$, and let $V$ be a local system of flat $\Lambda$-modules on $X$. The weight spectral sequence is a spectral sequence ($p \ge 0, q \ge 0$ are positive integers) - 

\begin{equation}  
E_{1}^{p, q} = \bigoplus_{i \ge max(0, -p)} H^{q-2i}(X_{\bar{s}}^{(p+2i)}, V(-i)) \Rightarrow H^{p+q}(X_{\bar{K}}, V)
\end{equation}
that computes the cohomology of the generic fiber in terms of the cohomology of the special fiber. This is equivariant for the action of $Gal(\bar{F_0} / F_0)$ on both sides, and the differentials commute with this action. Note that the groups $E^{p, q}_{1}$ vanish for $q < 0$ and $q > 2n$. The sequence can be seen to arise from the complex of nearby cycles $R\Psi V$. Let us briefly recall the construction in \cite{Saito}. Consider the following diagram -

\begin{equation}
\begin{array}{ccccc}

X_{\bar{s}}  & \xrightarrow{\bar{i}}  & X_{\OS_{\bar{F_0}}} & \xleftarrow{\bar{j}} & X_{\bar{\eta}} \\
\downarrow &  & \downarrow & & \downarrow \\
X_{s}  & \xrightarrow{i}  & X_{\OS_{F_0}} & \xleftarrow{j} & X_{\eta} \\

\end{array}
\end{equation}

The complex $R\Psi V := \bar{i}^{\ast}R\bar{j}_{\ast} V$ is the nearby cycles complex that receives an action of the inertia group $I_{F_0} \subset Gal(\bar{F_0} / F_0)$. Let $T$ be an element of $I_{F_0}$ that maps to a generator of $\Z_{\ell}(1)$ under the isomorphism $t_{\ell} : I_{F_0} \rightarrow \Z_{\ell}(1)$. Then $T-1$ induces a nilpotent endomorphism of $R\Psi V$, which we denote by $\nu$. Then we have the following (\cite[Section 2]{Saito}) -

\begin{prop}
\begin{enumerate}
    \item $R \Psi V$ lies in the abelian subcategory $\mathrm{Perv}(X_{\bar{s}}, \Lambda)[-n]$ of $-n$-shifted perverse sheaves with coefficients $\Lambda$.
    
    \item Let $M_{\bullet}$ denote the increasing monodromy filtration of $\nu$. For each positive integer $p$, let $a_p : X_{\bar{s}}^{(p)} \rightarrow X_{\bar{s}}$ be the canonical map. Then, for each integer $r \ge 0 $, there is a canonical isomorphism \begin{equation*}
        \bigoplus_{p-q = r} a_{p+q. \ast} V(-p)[-(p+q)] \cong \mathrm{Gr}^{M}_{r} R \Psi V.
    \end{equation*}
\end{enumerate}
\end{prop}

The weight spectral sequence is now the spectral sequence canonically associated to the filtered object $R \Psi V$. The first row of this spectral sequence can be computed using \cite[Proposition 2.10]{Saito}. We recall the result briefly below.

We recall the construction of a simplicial complex $\mathcal{K}$ whose simplicial cohomology computes the first row on the $E_1$-page of the spectral sequence. We have $$ E_{1}^{p, 0} = H^{0}(X_{\bar{s}}^{(p)}, V) \cong \bigoplus_{|A| = p+1} H^{0}(X_{A, \bar{s}}, V),$$ and unraveling the definition, the differential $$ d_{1}^{p, 0} : E_{1}^{p, 0} \rightarrow E_{1}^{p+1, 0} $$ can be seen to be the sum of the canonical pullbacks $ i_{A, A^{'}} : H^{0}(X_{A, \bar{s}}, V) \rightarrow H^{0}(X_{A^{'}, \bar{s}}, V)$ each multiplied by precisely the sign $\epsilon(A, A^{'})$ defined in (\ref{sign}). The simplicial complex $\mathcal{K}$ is defined as follows. The vertices of $\mathcal{K}$ are in bijection with the $X_i$, and the set $\{ X_{i_1}, \ldots, X_{i_r} \}$ corresponds to a simplex $\sigma_{A}$ if and only if for $A= \{i_1, \ldots, i_r \}, X_A$ is nonempty. The coefficient system $\mathcal{V}$ corresponding to $V$ is simply given by the global sections of the special fiber: $\sigma_{A} \rightarrow H^{0}(X_{A, \bar{s}}, V)$. We denote by $C^{\bullet}(\mathcal{K}, \mathcal{V})$ the complex computing the simplicial cohomology of $\mathcal{K}$ with coefficients in $\mathcal{V}$. Thus, by definition, $$C^{r}(\mathcal{K}, \mathcal{V}) = \bigoplus_{A \subset \{1, \ldots, h \} } H^{0}(X_{A, \bar{s}}, V) $$ where the sum is over subsets $A$ of cardinality $r+1$. The differential $d_{r} : C^{r}(\mathcal{K}, \mathcal{V}) \rightarrow C^{r+1}(\mathcal{K}, \mathcal{V})$ is given by the direct sum of restriction maps $$ res_{A, A^{'}} : H^{0}(X_{A, \bar{s}}, V) \rightarrow H^{0}(X_{A^{'}, \bar{s}}, V) $$ each multiplied by the sign $\epsilon(A, A^{'})$ as before, and we have the following result (follows from \cite[Proposition 2.10]{Saito})-

\begin{prop}
There is a canonical isomorphism of complexes $E^{\bullet, 0}_{1} \cong C^{\bullet}(\mathcal{K}, \mathcal{V})$.
\end{prop}

\section{Cohomology, degeneration, and level-raising}
\label{sec:cohomology}

In this section, we use the description of irreducible components of the special fiber of the varieties $Sh(G, U)$ given by the $p$-adic uniformization result of Corollary \ref{main}, to compute first rows on $E_1$ and $E_2$ pages of the weight spectral sequence in terms of the spaces of automorphic forms $\mathcal{A}(U^{\omega}, W_{\bm{\mu}, k})$. Then we prove a degeneration result following the use of a trick from \cite{Thorne} inspired by the use of weights in characteristic $0$ to show that the spectral sequence degenerates at $E_2$ level. Finally, we put everything together to obtain a level-raising result for $I$ and then $GL_n$.

\subsection{Automorphic forms and weight spectral sequence}
Recall that for $\SGU$ satisfying hypothesis \ref{hyp:mild}, we have the $p$-adic uniformization result from Corollary \ref{main} -

\begin{equation} \label{uni}
I(\Q) \backslash [\mathcal{M} \times G(\A^{\omega, \infty}) / U^{\omega} \cong (Sh(G, U) \otimes_{\OS_{R_{u}}} \OS_{\breve{R_{u}}})^{\wedge}.
\end{equation}
By the theory of Bruhat-Tits building $BT$ of $PGL_{n}(E_{\omega})$ and the Drinfeld upper half plane, the set of irreducible components in the special fiber of $\mathcal{M}$ is in bijection with the set $BT(0) \times \Z$, where $BT(i)$ for $i \in \mathbb{N}$ refers to the set of simplices of $BT$ of dimension $i$. $BT(0)$ corresponds to homothety classes of $\OS_{F_{\nu}}$-lattices $M \subset E_{\omega}^{n}$. We can in fact write down this set concretely in terms of the groups $U_i$ defined in section $5.2$. These are maximal compact subgroups that stabilize the distinct vertices, say $x_0, \ldots, x_{n-1}$, of the closure of the unique chamber of $BT$ fixed by $B$, and their intersection is precisely $B$. We can define a coloring map $\kappa : BT(0) \times \Z \rightarrow \Z / n \Z $ by sending $(M, b)$ to $\kappa(M, b) = log_{q}[M : \OS_{F}^{n}] + b$, where $q$ is the cardinality of the residue field. Then $\kappa(x_i, 0) = i$ induces an isomorphism of $GL_n(E_{\omega})$-sets $$ BT(0) \times \Z \cong \bigsqcup^{n-1}_{i = 0} GL_n(E_{\omega}) / U_i.$$ Furthermore, for each $i = 0, \ldots, n-1$, there is a bijection between the set of nonempty $(i+1)$-fold intersections of irreducible components of the special fiber of $\mathcal{M}$ and the set $$BT(i) \times \Z \cong \bigsqcup_{A \subset \{0, \ldots, n-1 \} } GL_n(E_{\omega}) / U_{A}.$$ The disjoint union runs over subsets $A$ of cardinality $i+1$. From this information, we see from (\ref{uni}), since $U$ is sufficiently small, that the irreducible components of the special fiber are in bijection with the set \begin{equation} \label{components} I(\Q) \setminus [BT(0) \times GL_n(E_{\omega}) / GL_n(E_{\omega})^{0} \times I(\A^{\omega, \infty}) / U^{\omega} ] \cong \bigsqcup_{i=0}^{n-1} I(\Q) \setminus I(\A^{\infty}) / U_{i}U^{\omega}.\end{equation} Thus, for each $i = 0, \ldots, n-1$, there is a bijection $$ \pi_{0}(Sh(G,U)_{\bar{s}}^{(i)}) \cong \bigsqcup_{A \subset \{ 0, \ldots, n-1 \} } I(\Q) \setminus I(\A^{\infty}) / U_{A}U^{\omega},$$ the union running over subsets $A$ of cardinality $i+1$. 

To write down the weight spectral sequence for $\SGU$, we choose a partial ordering on the set of irreducible components of the special fiber, i.e. on the $LHS$ of (\ref{components}), given by the pullback under $\kappa$ of the partial ordering $0 \le 1 \le \ldots \le n-1$ defined on $\Z / n \Z$. Let $E_{1}^{p, q} \Rightarrow H^{p+q}(\SGU_{\bar{\eta}}, V_{\bm{\mu}, k})$ be the weight spectral sequence for $\SGU$ and the local systems $V_{\bm{\mu}, k}$.

\begin{prop}
\label{autform}
\begin{enumerate}

\item For each $ i = 0, \ldots, n-1 $, there is a canonical isomorphism $$ E_{1}^{i, 0} \cong \bigoplus_{A \subset \{ 0, \ldots, n-1 \} } \mathcal{A}(U^{\omega}U_{A}, W_{\bm{\mu}, k})$$ where the direct sum runs over the set of all subsets $A$ of cardinality $i+1$.

\item There is a canonical isomorphism of complexes $$ E_{1}^{\bullet, 0} \cong C^{\bullet}(\mathcal{A}(U^{\omega}, W_{\bm{\mu}, k})) $$ and therefore for each $i = 0, \ldots, n-1$, $$E_{2}^{i, 0} \cong H^{i}(\mathcal{A}(U^{\omega}, W_{\bm{\mu}, k})).$$

\end{enumerate}
\end{prop}

\begin{proof}
 By definition, we have $ E_{1}^{i, 0} = H^{0}(Sh(G, U)^{(i)}_{\bar{s}}, V_{\bm{\mu}, k}) $. We have seen previously that this space of global sections can be identified with the space of functions $f : \pi_{0}(Sh(G, U(1))^{(i)}_{\bar{s}}) \rightarrow W_{\bm{\mu}, k}$ satisfying the transformation relation $f(C \sigma) = \sigma^{-1}f(C)$ for all $C \in  \pi_{0}(Sh(G, U(1))^{(i)}_{\bar{s}}), \sigma \in U$. The first isomorphism now follows from our identification of the set  $\pi_{0}(Sh(G, U(1))^{(i)}_{\bar{s}})$ with the disjoint union $\bigsqcup_{A \subset \{ 0, \ldots, n-1 \} } I(\Q) \setminus I(\A^{\infty}) / U_{A}U(1)^{\omega}$ and the definition of the spaces $\mathcal{A}(U^{\omega}U_A, W_{\bm{\mu}, k})$.

We note that under this isomorphism, the restriction maps of sections $H^{0}$ correspond to the natural inclusions of spaces of automorphic forms $\mathcal{A}(U^{\omega}U_A, W_{\bm{\mu}, k}) \rightarrow \mathcal{A}(U^{\omega}U_{A^{'}}, W_{\bm{\mu}, k})$, and the signs $\epsilon(A, A^{'})$ agree by their definition. This implies that the differentials of the two complexes in the second part correspond under the isomorphism of the first part, giving us the desired isomorphisms of complexes.
\end{proof}

\subsection{Degeneration}
We now describe the degeneration result that is at the heart of this method. To this end, we first need to define another descent of the scheme $Sh(G, U) \otimes_{\OS_{R_{u}}} \OS_{\breve{R_{u}}}$ on which the action of Frobenius in characteristic $\ell$ is given by scalar weights. We define $$ Sh(G, U)^{split} := I(\Q) \setminus [\mathcal{M}^{split}_{E_{\omega}} \times G(\A^{\omega, \infty}) / U^{\omega}]. $$ The local systems $V_{\bm{\mu}, \Lambda}$ where $\Lambda = K, \OS,$ or $\OS / \lambda^m$ admit descents to the scheme $Sh(G, U)^{split}$ by following the same definitions as before. $Sh(G, U)$ is not the same descent as $Sh(G, U)^{split}$ (in particular, Galois actions differ), but they become isomorphic after extension of scalars to $\breve{R_{u}}$. 
\begin{lem}
\begin{enumerate}
\item The pullback of $V_{\bm{\mu}, k}^{split}$ to any irreducible component of the special fiber of $Sh(G, U)^{split}$ is a constant sheaf.

\item If $Y_1, \ldots, Y_s$ are irreducible components of the special fiber, then the action of the Frobenius element is by the scalar $q_{\nu}^{i/2}$ on the group $H^{i}((Y_1 \cap \ldots \cap Y_s)_{\bar{s}}, V^{split}_{\bm{\mu}, k}).$
\end{enumerate}
\end{lem}

\begin{proof}
Let $Y \subset Sh(G, U(1))^{split}$ be an irreducible component of the special fiber. If $\pi : Sh(G, U(1))^{split} \rightarrow Sh(G, U)^{split}$ is the natural projection, then the restriction $\pi|_{Y}$ induces an isomorphism from $Y$ to its image in $Sh(G, U)^{split}$. Pullback of $V_{\bm{\mu}, k}$ under the inverse of this isomorphism gives the first part. The second part then follows from the first and Proposition $3.1$ from \cite{Thorne}.
\end{proof}

We recall that the group $H^{i}$ appearing in the second part above is $0$ if $i$ is odd.

\begin{prop}
\label{degenerate}
Let $r = 2s+1$. With notation and setup as above, the differentials $$ d_{r}^{p, q} : E_{r}^{p, q} \rightarrow E_{r}^{p+r, q+1-r} $$ are all zero as long as $q_{\nu}^{s} \not \equiv 1$ mod $\ell$.
\end{prop}

\begin{proof}
The weight spectral sequence of a pair $(X, V)$ only depends upon the extension of scalars $(X \otimes_{\OS_{F_{0}}} O_{\breve{F_0}}, V)$ as a spectral sequence of abelian groups (though  \textit{not} as a sequence of Galois modules). This implies that it is sufficient to prove the assertion for the pair $(Sh(G, U)^{split}, V_{\bm{\mu}, k}^{split})$. This then follows from the fact that the differentials in the weight spectral sequence are Galois equivariant, and the previous lemma.
\end{proof}

\begin{cor}
\label{injection}
Suppose that $\ell$ is a banal characteristic for $GL_n(E_{\omega})$. Then the weight spectral sequence for the pair $(Sh(G, U), V_{\bm{\mu}. k})$ degenerates at the $E_2$-level, and consequently there is an injection, equivariant for the prime-to-$p$ Hecke algebra : $$ H^{i}(\mathcal{A}(U^{\omega}, W_{\bm{\mu}, k})) \hookrightarrow H^{i}(Sh(G, U)_{\bar{R}_{u}}, V_{\bm{\mu}, k}).$$
\end{cor}

\begin{proof}
By the banal characteristic hypothesis, $\ell$ is coprime to the pro-order of $GL_n(E_{\omega})$. This implies that none of the integers $q_{\nu}, q_{\nu}^2, \ldots, q_{\nu}^{n-1}$ are congruent to $1$ modulo $\ell$. Since the cohomology of $Sh(G, U)$ can only non-trivially exist in the range $0, \ldots, 2(n-1)$, the groups $E_{r}^{p, q}$ can be nonzero only for that range. Then it follows from Proposition \ref{degenerate}, that all the differentials for $r \ge 2$ are zero, giving the first assertion. In particular, as noted before, since the weight spectral sequence is associated to the filtered complex $R\Psi V$, we have $E^{q, 0}_{\infty} \hookrightarrow H^{q}(Sh(G, U)_{\bar{R}_{u}}, V_{\bm{\mu}, k})$. The claim now follows from degeneration and Proposition \ref{autform}.
\end{proof}

\subsection{A level-raising congruence}
We briefly recall the level-raising congruence $(4.2)$ in \cite{Thorne} which will be applied in our setting in the proof of our main result. 

Let $L$ be a finite extension of $\Q_p$. Let $V$ be an admissible $GL_n(L)$-module over $\bar{\Q}_{\ell}$. We say that $V$ admits an integral structure if there exists a $GL_n(L)$-invariant $\bar{\Z}_{\ell}$-lattice $\Lambda \subset V$ such that $\Lambda \otimes_{\bar{\Z}_{\ell}} \bar{\Q}_{\ell} \cong V$. If $V$ admits an integral structure, then the reduction of $\Lambda$ modulo the maximal ideal of $\bar{\Z}_{\ell}$ gives an admissible representation of $GL_n(L)$ over $\bar{\mathbb{F}}_{\ell}$, which has finite length and its Jordan-Holder factors are independent of the chosen lattice $\Lambda$. (See \cite{Vigneras} for more information on integral structures.)

Let $T \subset P \subset GL_n(L)$ denote the standard maximal torus and standard Borel, with $Iw \subset GL_n(L)$ the Iwahori subgroup. Let $T_0 \subset T$ denote the unique maximal compact subgroup. Let $W$ denote the Weyl group. If $\rho$ is an admissible representation of $GL_n(L)$ over $\bar{\Q}_{\ell}$, then by the theory of Bernstein center, the Iwahori invariants $\rho^{Iw}$ of $\rho$ can be endowed with an action of the algebra $\bar{\Q}_{\ell}[T / T_0] \cong \bar{\Q}_{\ell}[X_1, X_1^{-1}, \ldots, X_n, X_{n}^{-1}]$. We write $t_i := e_{i}(X_1, \ldots, X_n) \in \bar{\Q}_{\ell}[T / T_0]^{W}$, where $e_i$ is the symmetric polynomial of degree $i$ in $n$ variables. These elements are fixed by the natural action of the Weyl group on $\bar{\Q}_{\ell}[T / T_0]$.

If $\rho$ is an irreducible admissible representation of $GL_n(L)$ over $\bar{\Q}_{\ell}$ with $\rho^{Iw} \not = 0$, then we say that $\rho$ satisfies the level-raising congruence if there exists $\bar{\alpha} \in \bar{\mathbb{F}}_{\ell}^{\times}$ such that, for each $i = 1, 2, \ldots, n$, the eigenvalue $\alpha_i$ of $t_i$ on $\rho^{Iw}$ satisfies the congruence \begin{equation} 
\label{eq:LRC}
\alpha_i \equiv \bar{\alpha} e_i(q^{\frac{n-1}{2}}, \ldots, q^{\frac{1-n}{2}}) \mod \mathfrak{m}_{\bar{\Z}_{\ell}}. \end{equation}

\subsection{A sign condition}

We define a sign condition here that will be a hypothesis on the automorphic representations that appear in the application. Let $\pi_1, \pi_2$ be conjugate self-dual cuspidal automorphic representations of $GL_{n_1}(\mathbb{A}_E), GL_{n_2}(\mathbb{A}_E)$ respectively, such that $ \pi = \pi_1 \boxplus \pi_2$ is regular algebraic. Then the representations $\pi_i |.|^{n_i - n / 2}$ are regular algebraic for $i = 1, 2$. In this case, define a tuple $\mathbf{a}^i = (a^{i}_{\tau})_{\tau \in Hom(E, \C)}$ by the requirement that $(a^{i}_{\tau, 1} + (n_i - n)/2, \ldots, (a^{i}_{\tau, n_i} + (n_i - n)/2)_{\tau \in Hom(E, \C)} $ equals the infinity type of $\pi_i |.|^{n_i - n / 2}$ for $i = 1, 2$. Then define $\mathbf{b} = (b_{\tau})_{\tau \in Hom(E, \C)}$ by the concatenation formula $$ (b_{\tau, 1}, \ldots, b_{\tau, 2}) = (a^{1}_{\tau, 1}, \ldots, a^{1}_{\tau, n_1}, a^{2}_{\tau, 1}, \ldots, a^{2}_{\tau, n_2}). $$ Let $S_n$ denote the permutation group on $\{ 1, 2, \ldots, n\}$. Then, there is a unique tuple $\mathbf{w} = (w_{\tau})_{\tau \in Hom(E, \C)} \in S_n^{Hom(E, \C)}$ such that, for each $\tau \in Hom(E, \C)$, the infinity type of $\pi$ is $(b_{\tau, w_{\tau}(1)}, \ldots, b_{\tau, w_{\tau}(n)})_{\tau \in Hom(E, \C)}$. We say that $\pi = \pi_1 \boxplus \pi_2$ satisfies the sign condition if the following condition holds true. Recall that we have $\tau'_{1}, \ldots, \tau'_{d}$ as embeddings at $\infty$ of $F$ and chosen $\tau_1, \ldots, \tau_d$ embeddings of $E$ extending them. Then we require that \begin{equation} \label{eq:sgn} \prod_{i = 1}^{d} sgn( w_{\tau_d}) = 1. \end{equation}

\subsection{Level-raising for unitary groups}
We now apply Corollary \ref{injection} to the problem of level-raising. For the rest of this subsection, we work in the notation from section $5$. We know that the contribution at $\omega$ to $I(\Q_p)$ is isomorphic to $GL_n(E_{\omega})$. We assume that under this isomorphism, $U_{\omega} \cong B$, where $B$, as before, is the standard Iwahori subgroup of $GL_n(E_{\omega})$. We denote by $U'_{\omega}$ the unique maximal compact subgroup at $\omega$. We then have the following level-raising theorem :

\begin{thm}[Level-raising for $I$]
\label{raise}
Suppose that $\sigma$ is as in section $5$, that is, If $\sigma$ is an automorphic representation of $I(\A)$ such that $(\sigma^{\infty})^{U} \neq 0$ and $\sigma_{\infty} \cong W^{\vee}_{\bm{\mu}, \C}$. Let $m_{\sigma}$ the associated maximal ideal of the Hecke algebra $\mathbb{T}_{T}^{univ}$. Suppose that the following hypotheses hold :
\begin{enumerate}
\item $U^{\omega}$ is a sufficiently small open compact subgroup of $I(\A^{\omega, \infty})$.

\item If $\sigma^{'}$ is another automorphic representation contributing to the space $\mathcal{A}(U, W_{\bm{\mu}, \mathscr{O}})_{m_{\sigma}}$, then $\sigma^{'}$ considered as a representation of $GL_n(E_{\omega})$ is a subquotient of a parabolic induction $n-Ind_{Q}^{G} St_{a}(\alpha) \otimes St_{b}(\beta)$ for some partition $a+b = n$ and $Q$ the standard parabolic corresponding to this partition.

\item $\iota^{-1} \sigma_{\omega}$ satisfies the level-raising congruence \ref{eq:LRC}.

\item $\bm{\mu}$ is $\ell$-small and $\ell$ is a banal characteristic for $GL_n(E_{\omega})$.

\item The groups $H^{n-2}(Sh(G, U'_{\omega}U^{\omega})_{\bar{R}_{u}}, V^{\vee}_{\bm{\mu}, k})$ and $H^{n-2}(Sh(G, U^{'}_{\omega}U^{\omega})_{\bar{R}_{u}}, V_{\bm{\mu}, k})$ are zero.
\end{enumerate}

Then there exists another irreducible constituent $\sigma^{'}$ contributing to the space $\mathcal{A}(U, W_{\bm{\mu}, \mathscr{O}})_{m_{\sigma}}$, such that $\sigma^{'}$ is an unramified twist of the Steinberg representation, i.e. $\sigma^{'}$ is obtained by raising the level of $\sigma$.
\end{thm}

\begin{proof}
By hypothesis $5$ and Corollary \ref{injection}, the groups, in characteristic $\ell$, $H^{n-2}(\mathcal{A}(U^{\omega}, W_{\bm{\mu}, k}))$ and $H^{n-2}(\mathcal{A}(U^{\omega}, W^{\vee}_{\bm{\mu}, k}))$ are zero. On the other hand, in characteristic $0$, we have a perfect $GL_n(E_{\omega})$-equivariant pairing $$ \langle ., . \rangle : \mathcal{A}(U^{\omega}, W_{\bm{\mu}, \mathscr{O}}) \times \mathcal{A}(U^{\omega}, W^{\vee}_{\bm{\mu}, \mathscr{O}}) \rightarrow \OS. $$ For $V \subset B$ an open compact subgroup, and $f_1 \in \mathcal{A}(U^{\omega}V, W_{\bm{\mu}, \mathscr{O}}), f_2 \in \mathcal{A}(U^{\omega}V, W^{\vee}_{\bm{\mu}, \mathscr{O}})$, we define $$ \langle f_1, f_2 \rangle = \frac{1}{[B:V]} \sum_{x \in I(\Q) \backslash I(\A^{\infty}) / U^{\omega}V } (f_1(x), f_2(x)). $$ This is independent of the choice of $V$ and equivariant for the action of $GL_n(E_{\omega})$. The action of Hecke algebra on $\mathcal{A}(U^{\omega}, W_{\bm{\mu}, \mathscr{O}})$ gives a decomposition $$ \mathcal{A}(U^{\omega}, W_{\bm{\mu}, \mathscr{O}}) = \mathcal{A}(U^{\omega}, W_{\bm{\mu}, \mathscr{O}})_{m_{\sigma}} \oplus C $$ for some $C$. The result then follows, appealing to the level-raising formalism developed in \cite[Section 4]{Thorne} with $M = \mathcal{A}(U^{\omega}, W_{\bm{\mu}, \mathscr{O}})_{m_{\sigma}}$ and $N$ being the annihilator of $C$ under the pairing $\langle .,. \rangle$.
\end{proof}

\section{Application to $\GLn$}
In this section, we deduce our main theorem. $E$ is a CM field as always with maximal totally real subfield $F$, $[F:\Q] = d$ is even. The level-raising prime $p$ is assumed to be unramified in $F$, with a place $\nu$ above it of $F$ assumed to be split in $E$ as $\nu = \omega . \omega^{c}$. Let $n \ge 3$ be an integer, and $\ell \not = p$ a prime. We fix an isomorphism $\iota : \bar{\Q}_{\ell} \cong \C $ as before.

Let $n_1, n_2$ be positive integers with $ n = n_1 + n_2 $. Suppose that $\pi_1, \pi_2$ are conjugate self-dual cuspidal automorphic representations of $GL_{n_1}(\A_E)$ and $GL_{n_2}(\A_E)$ respectively such that $\pi = \pi_1 \boxplus \pi_2$ is regular algebraic. By work of many people (\cite[Theorem 1.7]{CH} and \cite[Theorem 1.1]{Caraiani} for our particular case) there is associated to $\pi$ a continuous semisimple Galois representation $r_{\iota}(\pi) : G_{E} \rightarrow GL_n(\bar{\Q}_{\ell})$.

\begin{thm}
\label{final}
With $\pi$ as above, suppose that $\iota^{-1} \pi_{\omega}$ satisfies the level-raising congruence \ref{eq:LRC}. Suppose further that the following hypotheses hold.

\begin{enumerate}
\item If $t_{\ell}$ is a generator of the $\ell$-th part of the tame inertia group at $\omega$, then $\overline{r_{\iota}(\pi)}(t_{\ell})$ is a unipotent matrix with exactly two Jordan blocks.

\item $\ell$ is a banal characteristic for $GL_{n}(E_{\omega})$.

\item The weight $\bm{\lambda} = (\lambda_{\tau})_{\tau : E \hookrightarrow \C}$ of $\pi$ satisfy the following properties :
          
            \begin{itemize}

            \item For each $\tau$, and for each $0 \le i < j \le n$, we have $0 < \lambda_{\tau, i} - \lambda_{\tau, j} < \ell$.

            \item There exists an isomorphism $\iota_{p} : \overline{\Q}_p \cong \C$ such that the following inequalities hold : 

		$$ 2n + \sum_{\tau : E \hookrightarrow \C} \sum_{j = 1}^{n} (\lambda_{\tau, j} - 2 \left \lfloor{ \lambda_{\tau, n}/2 }\right \rfloor ) \le \ell,$$
		$$ 2n + \sum_{\tau : E \hookrightarrow \C} \sum_{j = 1}^{n} (2 \left \lceil{ \lambda_{\tau, 1}/2 }\right \rceil  - \lambda_{\tau, n+1-j}) \le \ell.$$

            \end{itemize}

\item If $\pi$ is ramified at a place $\gamma$ of $E$, then $\gamma$ is split over $F$.

\item $\pi$ is unramified at the primes of $E$ lying above $\ell$, and $\ell$ is unramified in $E$.

\item $\pi = \pi_1 \boxplus \pi_2$ satisfies the sign condition \ref{eq:sgn}, $n_1 \not = n_2$, and $n_1n_2$ is even.

\end{enumerate}
Then there exists a regular algebraic conjugate self-dual cuspidal automorphic representation $\Pi$ of $GL_n(\A_E)$ of weight $\bm{\lambda}$ such that there is an associated Galois representation $r_{\iota}(\Pi)$ to it satisfying $\overline{r_{\iota}(\pi)} \cong \overline{r_{\iota}(\Pi)}$ and $\Pi_{\omega}$ is an unramified twist of the Steinberg representation. If the places of $F$ above $\ell$ are split in $E$, and $\pi$ is $\iota$-ordinary in the sense of \cite[Definition 5.1.2]{Geraghty}, then we can even assume that $\Pi$ is also $\iota$-ordinary.
\end{thm}

\begin{proof}
This follows from Theorem \ref{raise}, once we make all the right identifications. By \cite[Proposition 2.9]{CT}, and \cite[Corollaire 5.3]{Labesse}, we see that there exists an automorphic representation $\sigma$ of $I$ such that $\pi$ is the base change of $\sigma$. Once we have $\sigma$ in our hands, it is a matter of checking that various hypotheses of Theorem \ref{raise} are satisfied. Only the hypotheses $(2)$ and $(5)$ there need checking. Hypothesis $(5)$ of that theorem follows from the torsion vanishing result of Lan and Suh \cite[Theorem 8.12]{LS} combined with the inequalities in hypothesis $(3)$ of this theorem, and Proposition \ref{prop:compo} that shows the fact that $S(G,U)$ is a union of connected components of unitary similitude Shimura varieties.  One could also appeal to the torsion vanishing result of \cite{Shin} after putting some extra hypotheses. The hypothesis $(2)$ in that theorem is checked in a similar way as that of Theorem $7.1$ of \cite{Thorne}, from our hypothesis $(1)$ on Jordan blocks. If $\sigma^{'}$ is the representation whose existence is guaranteed by Theorem \ref{raise} for $\sigma$ as above, we use the reverse direction of the aforementioned base change results to obtain a representation $\Pi$ having required properties. The last assertion about ordinary property follows after enlarging the Hecke algebra $\mathbb{T}_{T}^{univ}$ to contain the $U_{\ell}$ operators.
\end{proof}

\begin{rem}
The analogue of Theorem $1.1$ from \cite{Thorne} with the hypothesis on $p$ being `$p$ is unramified in $F$' as opposed to `$p$ is inert in $F$' follows from the above theorem in the same way as the proof there. Note that $n_1$ and $n_2$ are assumed to be distinct even integers for this theorem.
\end{rem} 

\begin{rem}
Caraiani and Scholze \cite{CS} have proven stronger torsion vanishing results than \cite{LS} or \cite{Shin} for cohomology of Shimura varieties. Our theorem could be stated in terms of those results, too. Although we note that they assume that the CM field $E$ contains an imaginary quadratic field.
\end{rem}

%    Bibliographies can be prepared with BibTeX using amsplain,
%    amsalpha, or (for "historical" overviews) natbib style.
\bibliographystyle{amsplain}

\begin{thebibliography}{CHT08}

\bibitem[AC89]{AC}
J. Arthur and L. Clozel,
\emph{Simple Algebras, Base Change, and the Advanced Theory of the Trace Formula,}
Annals of Mathematics Studies $\mathbf{120}$ (1989).

\bibitem[AT19]{AT}
C. Anastassiades and J. A. Thorne,
\emph{Raising the level of automorphic representations of $GL_{2n}$ of Unitary type,}
Preprint, available at https://arxiv.org/abs/1912.11267 .

\bibitem[Ca12]{Caraiani}
A.~Caraiani,
\emph{Local–global compatibility and the action of monodromy on nearby
cycles,}
Duke Math. J. 161 (12), 2311 -– 2413, (2012).

\bibitem[Ca86]{Carayol} 
H.~Carayol,
\emph{Sur la mauvaise r\'{e}duction des courbes de {Shimura},}
Compositio. Math. 59(2), 151 -- 230, (1986).

\bibitem[Ce76]{Cerednik} 
I.~\v{C}erednik,
\emph{Towers of algebraic curves that can be uniformized by discrete subgroups of $PGL_2(K_{\omega}) \times E$,}
Math. USSR-Sb. 28 (1976), no. 2, 187 -- 215 (1978).

\bibitem[CH13]{CH}
G. Chenevier and M. Harris, 
\emph{Construction of automorphic Galois representations, II,} 
Cambridge J. Math., 1, 53 -– 73 (2013).

\bibitem[Cl91]{Clozel}
L. Clozel,
\emph{Représentations Galoisiennes associées aux représentations automorphes autoduales de $GL(n)$,}
Publ. Math. Inst. Hautes \'{E}tudes Sci. (73), 97 - 145 (1991).

\bibitem[CHT08]{CHT} 
L. Clozel, M. Harris, R. Taylor,
\emph{Automorphy for some $l$-adic lifts of automorphic mod $l$ Galois representations,}
Publ. Math. Inst. Hautes \'{E}tudes Sci. (108), 1 -- 181, (2008).
With Appendix A, summarizing unpublished work of Russ Mann, and Appendix B by Marie-France Vign\'{e}ras.

\bibitem[CS17]{CS}
A. Caraiani, P. Scholze,
\emph{On the generic part of the cohomology of compact unitary Shimura varieties,}
Annals of Mathematics (2) \textbf{186}, no. 3, 649 -- 766, (2017). 

\bibitem[CT14]{CT} L. Clozel, J. A. Thorne,
\emph{Level-raising and symmetric power functoriality, I,}
Compos. Math, 150(05), 729 -- 748, (2014).

\bibitem[De71]{Deligne} P. Deligne,
\emph{Travaux de Shimura,}
in S\'{e}minaire Bourbaki, 23\`{e}me ann\'{e}e (1970/71), Exp. No. 389', Springer, Berlin, Lecture Notes in Math., Vol. 244, 123 -- 165.

\bibitem[Dr74]{Drinfeld} V. Drinfeld,
\emph{Elliptic modules,}
Mat. Sb. (N.S.), 94(136):4(8), 594 -- 627, (1974).

\bibitem[Dr76]{Drinfeld02}
V. Drinfeld,
\emph{Coverings of $p$-adic symmetric regions,}
Functional Anal. Appl. $\mathbf{10}$, 107 -- 115 (1976).

\bibitem[DT94]{DT}
F. Diamond and R. Taylor,
\emph{Nonoptimal levels of mod $l$ modular representations,}
Invent. Math.115, no. 3, 435 –- 462, (1994).

\bibitem[Ge10]{Geraghty} D. Geraghty,
\emph{Modularity lifting theorems for ordinary Galois representations,}
PhD Thesis, Harvard University (2010).

\bibitem[Gr99]{Gross} B. Gross,
\emph{Algebraic modular forms,}
Israel J. Math. 113, 61 -- 93, (1999).

\bibitem[Ha13]{Harris} M. Harris,
\emph{The Taylor-Wiles method for coherent cohomology,}
J. Reine Angew. Math 679, 125 -- 153, (2013).

\bibitem[HT01]{HT} M. Harris and R. Taylor,
\emph{The geometry and cohomology of some simple Shimura varieties,}
volume 151 of Annals of Mathematics Studies, Princeton University Press, Princeton, NJ, 2001. With an appendix by Vladimir G. Berkovich.

\bibitem[Ko92A]{Kottwitz01} R. Kottwitz,
\emph{On the $\lambda$-adic representations associated to some simple Shimura varieties,}
Invent. Math. 108, 653 -- 665, (1992).

\bibitem[Ko92B]{Kottwitz02} R. Kottwitz,
\emph{Points on some Shimura varieties over finite fields,}
Jour. of the AMS 5, 373 -- 444, (1992).

\bibitem[La11]{Labesse} J.-P. Labesse,
\emph{Changement de base CM et s\'{e}ries discr\`{e}tes,}
in : Stab Trace Formula Shimura Var Arith Appl. Vol. 1 (Int Press, Somerville, MA, 2011), 429 -- 470.

\bibitem[LS12]{LS} K-W. Lan and J. Suh,
\emph{Vanishing theorems for torsion automorphic sheaves on compact PEL-type Shimura varieties,}
Duke Math. J. 161(6), 1113 -- 1170, (2012).

\bibitem[Mi05]{Milne} J. Milne,
\emph{Introduction to Shimura Varieties,}
Harmonic analysis, the trace formula, and Shimura varieties, volume 4 of Clay Math. Proc. Amer. Math. Soc., Providence, RI, 2005.

\bibitem[NT19]{NT} J. Newton and J. A. Thorne,
\emph{Symmetric power functoriality for holomorphic modular forms,}
preprint, available at https://arxiv.org/abs/1912.11261 .

\bibitem[RZ96]{RZ} M. Rapoport and T. Zink,
\emph{Period spaces for $p$-divisible groups,}
volume 141 of Annals of Mathematics Studies, Princeton University Press, Princeton, NJ, 1996.

\bibitem[Ri83]{Ribet} K. A. Ribet,
\emph{Congruence relations between modular forms,}
Proceedings of the International Congress of Mathematicians, vol. 1, 2 (Warsaw, 1983) (PWN, Warsaw, 1984), 503 -- 514.

\bibitem[Sa03]{Saito} T. Saito,
\emph{Weight spectral sequences and independence of $l$,}
J. Inst. Math. Jussieu {\textbf 2} (4), 583 -- 634, (2003).

\bibitem[Sh15]{Shin} S. W. Shin,
\emph{Supercuspidal part of the mod $l$ cohomology of $GU(1,n-1)$-Shimura varieties,}
J. Reine Angew. Math. 705, 1 -- 21, (2015). 

\bibitem[Th14]{Thorne} J. Thorne,
\emph{Raising the level for $GL_n$,}
Forum of Mathematics $\Sigma$, Vol. 2, e16, (2014).

\bibitem[Va98A]{Varshavsky01} Y. Varshavsky,
\emph{$p$-adic uniformization of unitary Shimura varieties,} 
Inst. Hautes \'{E}tudes Sci. Publ. Math. 87, 57 -- 119, (1998).

\bibitem[Va98B]{Varshavsky02} Y. Varshavsky,
\emph{$p$-adic uniformization of unitary Shimura varieties. II.,} 
J. Differential Geom. 49, no. 1, 75 -- 113, (1998).

\bibitem[Vi94]{Vigneras} M.-F. Vign\'{e}ras,
\emph{Banal characteristic for reductive $p$-adic groups,}
J. Number Theory 47(3), 378 -- 397, (1994).


\bibitem[Zh14]{Zhang} W. Zhang,
\emph{Selmer groups and the indivisibility of Heegner points,}
Cambridge Journal of Mathematics, Vol. 2, No. 2, 191 -- 253, (2014).
\end{thebibliography}
%    Insert the bibliography data here.

\end{document}